\documentclass[12pt]{amsart}
\usepackage{amssymb,amsmath,amsthm,mathrsfs}
\usepackage[dvipsnames]{xcolor}
\usepackage{graphicx}
\usepackage{tabularx}
\usepackage{multirow}
\usepackage{makecell}
\usepackage{setspace}
\usepackage{float}
\usepackage{enumitem}
\usepackage{tikz-cd}
\usepackage[linktocpage=true,backref=page]{hyperref}
\hypersetup{
	colorlinks=true, 							%set true if you want colored links
	linktoc=all,     								%set to all if you want both sections and subsections linked
	linkcolor=Fuchsia,  							%choose some color if you want links to stand out
	citecolor=Emerald,
}
\calclayout
\allowdisplaybreaks
\restylefloat{table}
\newtheorem{theorem}{Theorem}[section]
\numberwithin{theorem}{section} \numberwithin{equation}{section}

\newtheorem{definition}[theorem]{Definition}

\theoremstyle{remark}
\newtheorem{remark}[theorem]{Remark}

\newcommand\scalemath[2]{\scalebox{#1}{\mbox{\ensuremath{\displaystyle #2}}}}
\title[Twists, Eisenstein series, and Instantons in Local Mirror Symmetry]{Twists, Eisenstein series, and Instantons in Local Mirror Symmetry}

\thanks{A.M. acknowledges support from the Simons Foundation through grant no.~202367.}
\author{Andreas Malmendier}
\address{Dept.\!~of Mathematics \& Statistics, Utah State University, Logan, UT 84322}
\email{andreas.malmendier@usu.edu}

\author{Michael T.  Schultz}
\address{Dept.\!~of Mathematics, Virginia Tech, Blacksburg, VA 24060} 
\email{michaelschultz@vt.edu}
\begin{document}
	\begin{abstract}
		We propose a mechanism for computing the genus zero invariants of local Calabi-Yau fourfolds arising as the total space of the canonical bundle of a ~rank-1 Fano threefold. Our method relies heavily on modular parameterizations of the associated Landau-Ginzburg model, as well as extension regulator classes and higher normal functions studied by Doran and Kerr in this setting. The generating function of the local invariants is proposed as a functional inverse of a weight-4 modular form of Eisenstein type that is naturally computed from the mirror data, in analogy with some known results for local threefolds that computes genus zero Gopakumar-Vafa invariants. Using the twist construction of Doran and the first named author, these two constructions are connected via Doran's generalized functional invariant map on the periods.
	\end{abstract}
	\keywords{Elliptic curves, K3 surfaces, Fano varieties, Picard-Fuchs equations, modular forms, local mirror symmetry}
	\subjclass[2020]{11F11, 14J27, 14J33, 33C60, 14N35}
	\maketitle
	\tableofcontents
	\section{Introduction}
	\label{s-intro}
	\par Modularity plays a subtle but striking role in mirror symmetry for Calabi-Yau varieties in low dimensions (less than three in classical mirror symmetry, and less than five in local mirror symmetry respectively). In classical mirror symmetry, this is due to complex structure moduli spaces of elliptic curves and K3 surfaces having natural and well-known connections to modular varieties; for the exact same reason, modularity naturally makes an appearance in local mirror symmetry in dimensions three and four. In contrast, the classical mirror symmetry computation \cite{MR1115626} for a smooth quintic threefold in $\mathbb{P}^4$ of the Yukawa coupling - a generating function whose Lambert expansion  coefficients are instanton numbers, or genus zero Gopakumar-Vafa (GV) invariants of the quintic threefold - appears naturally as a $q$-series derived from the mirror map, and is reminiscent of a weight-4 modular form of Eisenstein type. Despite this resemblance, it resists a direct interpretation in terms of classical modular objects \cite{movasati2024}. However, local mirror symmetry in dimension three exhibits a clearer modular structure. In certain cases, the moduli spaces of local Calabi-Yau threefolds are modular curves: the periods of the compact part of the mirror (an elliptic curve) are modular forms, the (closed) mirror map is a Hauptmodul, and the instanton generating functions are quasi-modular. Doran \cite{MR1779161} gave precise conditions in which the (closed) mirror map is a Hauptmodul. Other examples of modularity appearing in local symmetry include  \cite{MR3780269,MR2851155,MR1486340,MR1468700,zhou2014,MR3822918}, for example. This modularity reflects both the rigidity of low-dimensional Calabi-Yau moduli spaces and the possibility of extending local analytic computations near the large complex structure limit to a Zariski-dense open subset. It also has consequences for the asymptotic growth of Gromov-Witten invariants of local threefolds ~\cite[Cor. 5.3]{MR2851155}.
	\par The aim of this article is to extend this structure to local Calabi-Yau fourfolds, specifically those arising as the total space of the canonical bundle of a rank-1 Fano threefold. Our starting point is the work of Stienstra \cite{MR2282958}, who observed experimentally that a suitably normalized functional inverse of a natural weight-3 modular form reproduces instanton counts for local threefolds. A mathematical explanation of this phenomenon was subsequently provided by Doran and Kerr \cite{MR2851155}, and by Bloch, Kerr, and Vanhove \cite{MR3780269}, using higher normal functions and extension regulator classes. In the present setting, the mirror family of anticanonical K3 surfaces in a rank-1 Fano threefold forms a pencil of Picard rank-19 K3 surfaces with a canonical $\mathrm{M}_N$-polarization, and is thus a Shioda–Inose partner of a Kummer surface associated with a pair of $N$-isogenous elliptic curves \cite{MR1420220,MR2306141}. The lattice $\mathrm{M}_N$ is given by
	\begin{equation*}
		\mathrm{M}_N := \mathrm{H} \oplus \mathrm{E}_8(-1)^{\oplus 2} \oplus \langle -2N \rangle \, .
	\end{equation*}
	By results of Doran \cite[Thm. 5.13]{MR1738862}, \cite[Thm. 6]{MR1779161}, the modularity of these K3 families follows from that of the associated elliptic curve families. Combining this with the framework of Doran and Kerr, we extend Stienstra’s construction to the fourfold setting and propose an analogous mechanism for computing local invariants.
	\par The structure of this article is as follows. In \S \ref{s-MUM_Operators}, we review aspects of classical mirror symmetry and local mirror symmetry. In \S \ref{ss-local_threefolds}, we review known results related to modularity and enumerative invariants of local threefolds. Crucial for us is the work of Doran and Kerr \cite[\S 9, \S 10]{MR2851155}, who show that the modular forms that we discuss are always of Eisenstein type (Theorem \ref{thm-eisenstein}). Their work on higher normal functions and extension regulator classes allows for a natural relation between the closed mirror map and the local mirror map that we rely heavily on. For a Picard-Fuchs operator $\mathcal{L}_n$, we define the so-called \emph{quantum} operator
	\begin{equation*}
		\mathcal{D}_{n+1} = \theta \circ \mathcal{L}_n \, ,
	\end{equation*}
	where $\theta$ is the usual first order logarithmic derivative operator. The work of Doran and Kerr shows that this operator is ``dual'' in an appropriate sense to the usual operator $\mathcal{L}_n \circ \theta$ that appears in local mirror symmetry. The quantum operator is an important ingredient for computing the weight-4 Eisenstein series that are the analogue of those appearing in the local threefold story.
	\par In \S \ref{s-fano_mirror_symmetry} we describe results related to the modularity of Fano surfaces and threefolds. This is part of a larger story discussed recently by Doran et al \cite{MR4913468} that addresses the story for higher Picard rank than under consideration in this article. The mirror elliptic curves and lattice polarized K3 surfaces we consider can be linked by the so-called twist construction of Doran and the first named author \cite{MR4069107}, and studied further by the present authors \cite{MR4494119}. The utility of the twist construction is that it allows us to derive \emph{all} data from the Weierstrass model associated to the mirror elliptic curves. A similar story holds for other Fano threefolds that do not share this connection, but for the sake of brevity we do not list the relevant data as much of it appears elsewhere in the literature. In Appendices \ref{s-modular_DEs} and \ref{s-twist}, we review relevant ingredients of the modular parameterization of Picard-Fuchs equations and the twist construction. In Appendix \ref{ss-modular_data} we compile tables of modular data, Picard-Fuchs operators, and quantum operators that can be linked by the twist construction.
	\par Lastly, in \S \ref{s-genus_zero_invariants}, we provide our proposal for local fourfold invariants as a natural functional inverse of the weight-4 Eisenstein series associated with the mirror K3 surface (Definition \ref{def-weight_4_inverse_expansion}). In \S  \ref{sss-local_P3}, we demonstrate that we can recover known local fourfold invariants for local $\mathbb{P}^3$ (Equation \ref{eq-KP3_GW0_invariants}). Our computation is distinct from other known examples in the literature, such as by Klemm and Pandharipande \cite{MR2415462}. Our computation is linked, via the twist construction, to the analogous computation for local $\mathbb{P}^2$ and other examples studied by physicists, which we discuss in \S \ref{sss-local_P2}. From these experimental results, we conjecture in \S \ref{ss-conj} that similar results hold for all other families of rank-1 Fano threefolds. We compile in Table \ref{table-local_instantons} the coefficients of the leading order ``local instanton'' expansions associated to the Fano threefold families under consideration in the article, similar data is readily computable for the remaining families.
	\par In the context of this article, a Calabi-Yau (CY) variety is an $m$-dimensional complex algebraic variety $V$ with trivial canonical bundle, $K_V \cong \mathcal{O}_V$. When $m > 1$, we also require that $h^i(V,\mathcal{O}_V)=0$ for $i = 1,\dots, m-1$. Such a variety may be projective or open. A Fano variety $F$ is a smooth complex projective variety with ample anticanonical class. Given $F$ of dimension $n$, there are always two important CY varieties associated to $F$: anticanonical CYs $\check{V} \in |-K_F|$ of dimension $n-1$, and the open CY variety $K_F$ of dimension $n+1$, the total space of the canonical line bundle of $F$. We will consider both Fano surfaces (i.e., del Pezzo surfaces), and Fano threefolds of Picard rank $\rho(F)=1$. The del Pezzo surfaces $B_m$ are the surfaces with ample anticanonical class $-K$ of degree $(K,K) = m = 1,\dots,9$, while the relevant Fano threefolds are generic members of the 17 deformation classes of algebraic threefolds with ample anticanonical class such that the Picard lattice has rank-1. Such varieties were classified by Iskovskikh \cite{MR0463151,MR0503430}. We will not need a precise geometric description of them and defer the interested reader to Iskovskikh's work; see also Golyshev \cite[\S 6.3]{MR2306141} or Coates et al \cite{MR3470714}. Noether-Lefschetz type results imply that $\rho(F)=\rho(\check{V})=\rho(K_F)$ for generic $\check{V} \in |-K_F|$. Of course, the simplest such examples are $B_9=\mathbb{P}^2$ and $\mathbb{P}^3$, respectively, which contain the anticanonical cubic curves and quartic surfaces. 
	\medskip
	\begin{flushleft}
		\emph{Acknowledgment.} After an initial draft of this paper was posted, we received valuable feedback and insight from Fenglong You. We are grateful for his input and references to the literature \cite{MR3997137,MR3948687,MR4716739}. We also thank the referee for their many helpful comments, suggestions, and corrections, particularly for their bringing our attention to the work of Doran and Kerr \cite{MR2851155}. This work was partially supported by a grant from the Simons Foundation.
	\end{flushleft}
	\section{Geometry of maximally unipotent operators}
	\label{s-MUM_Operators}
	\par In this section we establish notation and recall relevant results, beginning with mirror symmetry for CY varieties in \S \ref{ss-classical_mirror_symmetry} and local mirror symmetry in \S \ref{ss-local_mirror_symmetry}. The mirror varieties under consideration in this article have geometry that is strongly governed by modular curves, and so we conduct a review in Appendix \ref{s-modular_DEs} of the classical Gauss-Schwarz theory and canonical elliptic fibrations defined over a modular curve. Known results relevant to modularity and mirror symmetry of local threefolds are stated in \S \ref{ss-local_threefolds}.
	\subsection{Classical mirror symmetry}
	\label{ss-classical_mirror_symmetry}
	\par Enumerative geometry of rational curves on a variety may be packaged into the geometric machinery of quantum cohomology. This turns the complex cohomology ring into a family of Frobenius algebras, each of whose product structure is known as a quantum product, and is determined by choosing a fixed cohomology class which then deforms the usual cup product, and encodes curve counts in the form of Gromov-Witten (GW) invariants. See, e.g., Manin \cite{MR1702284}. These numbers are generically rational, whose denominators are often large and reflect the reality that the moduli space of stable maps of rational curves with prescribed number of marked points and degree to a given variety is a Deligne-Mumford stack. Dubrovin \cite{MR1397274} showed that there is a pencil of canonically flat connections on a vector bundle over the cohomology space of the variety that governs the Frobenius structure; it is called the quantum $D$-module of the variety and has precisely one regular singular point and one irregular singular point. The small quantum $D$-module is the one obtained by restricting to the quantum $D$-module to the divisorial directions.
	\par For a CY variety $\check{V}$ of dimension three, the GW invariants can be combined in a systematic way to recover integer valued invariants that count the ``correct" number of curves of prescribed genus and degree on $\check{V}$ \cite{MR1204770,gv1998-1,gv1998-2}. Physicists call these instantons or BPS counts, while mathematicians call these Gopakumar-Vafa (GV) invariants. A generalization of these invariants for CY fourfolds was provided by Klemm \& Pandharipande \cite{MR2415462}. We will uniformly call such invariants instantons, or genus zero invariants through the remainder of this article. Mirror symmetry studies pairs $(\check{V},V)$ such that the quantum cohomology of $\check{V}$ can be effectively encoded in the classical geometry of $V$, such as the variation of Hodge structure (VHS). The latter is associated to the canonically flat Gauss-Manin connection \cite{MR229641,MR233825,MR282990}, which was shown by Griffiths and Deligne \cite{MR417174} using Hironaka's famous resolution of singularities \cite{MR199184} to possess only regular singular points. A brief summary of these ideas relevant to CY varieties appears in \S \ref{sss-VHS}. Thus, mirror symmetry provides in part a regularization scheme of the small quantum D-module. These are respectively called the A-side and B-side of mirror symmetry, in reference to Type IIA and Type IIB string theory.
	\par The classical mirror $V$ of a generic anticanonical hypersurface $\check{V} \subset \mathbb{P}^n=\mathbb{P}(X_0,\dots,X_n)$ of Picard rank 1 is well known to be determined by the family $\hat{V}$ of deformed Fermat hypersurfaces
	\begin{equation}
		\label{eq-deformed_fermat}
		\hat{V} : \quad	X_0^{n+1}+X_1^{n+1}+\cdots+X_n^{n+1}+(n+1) \lambda X_0 X_1 \cdots X_n=0 ,
	\end{equation}
	parameterized by $\lambda \in \mathbb{P}^1 - \{0,\infty\}$ via the Greene-Plesser orbifolding construction \cite{MR1059831}. Then (\ref{eq-deformed_fermat}) for $n-1=1$ is the Hesse pencil of elliptic curves, $n-1=2$ is the Dwork family of K3 surfaces of Picard rank-19, and $n-1=3$ is the three-dimensional CY hypersurface that arose in connection with the famous mirror quintic family studied by Candelas, de la Ossa, Green, and Parkes \cite{MR1115626}. As a reminder, the mirrors $V$ themselves are obtained from $\hat{V}$ via the orbifolding construction as follows. The finite abelian group $G_{n-1} := (\mathbb{Z}/(n+1)\mathbb{Z})^{n-1}$ acts via the generators given by roots of unity as $(X_0,X_j) \mapsto (\zeta_{n+1}^nX_0,\zeta_{n+1}X_j)$ for $j=1,\dots,n$ and $\zeta_{n+1} = \exp(\frac{2\pi i}{n+1})$. Then Equation (\ref{eq-deformed_fermat}) defining $\hat{V}$ is invariant under $G_{n-1}$, and the natural invariant coordinates 
	\begin{equation*}
		t := \frac{(-1)^{n+1}}{\lambda^{n+1}} \, , \quad  \, x_j := \frac{X_j^n}{(n+1)\lambda X_0 \cdots X_n} \, , \quad j= 1,\dots,n
	\end{equation*}
	descend to the quotient $\hat{V}/G_{n-1}$. The family $V$ parameterized by $t \in \mathbb{P}^1 - \{0,\dots,\infty\}$ is determined by the hypersurface in $\mathbb{P}^n=\mathbb{P}(x_0,\dots,x_n)$ from the natural (affine) relation of degree $n+1$
	\begin{equation}
		\label{eq-mirror_hypersurface}
		V:  \quad	x_1\cdots x_n(x_1+\cdots+x_n+1) + \frac{(-1)^{n+1}t}{(n+1)^{n+1}} = 0 \, .
	\end{equation}
	Batyrev and Borisov \cite{MR1416334} proved that the pencil $V$ is mirror to a generic hypersurface $\check{V} \in |\mathcal{O}(n+1)|$, in the sense that $h^{p,q}(\check{V})=h^{n-1-p,q}(V)$. Moreover, there is an analytic isomorphism called the mirror map of a neighborhood of the point $i\infty$ in the compactification of the complexified symplectic moduli space $\check{X}$ of $\check{V}$ with the point $t=0$ of the classical complex structure moduli space $X$ of $V$ of equal dimension $\dim(\check{X})=h^{1,1}(\check{V})=h^{n-1,1}(V)=\dim(X)$. The points are respectively called the large radius limit and large complex structure limit.
	\par The case $n-1=2$ is meaningful in the case of algebraic K3 surfaces, and was formally studied by Dolgachev \cite{MR1420220}, building off previous works of Nikulin and Pinkham, e.g., \cite{MR525944,MR429876}. There, it was shown that the N\'eron-Severi and transcendental lattices of mirror K3s satisfy $\mathrm{T}(V) \cong (\mathrm{Pic}(\check{V})\oplus\mathrm{H})^\perp$, where $\mathrm{H}$ is the canonical hyperbolic lattice of rank 2. Here the orthogonal complement is taken with respect to the ambient K3 lattice $\Lambda_{\mathrm{K}3} := \mathrm{H}^{\oplus 3} \oplus \mathrm{E}_8(-1)^{\oplus 2}$ of rank $22$ and signature $(3,19)$, with $\mathrm{E}_8(-1)$ is the negative definite root lattice associated to that of the same name. In the case of a generic quartic K3 surface $\check{V} \subset \mathbb{P}^3$, we have $\mathrm{Pic}(\check{V}) = \langle 4 \rangle$, and it follows from the lattice relation above that the mirror $V$ is polarized by the lattice $\mathrm{M}_2 := \mathrm{H} \oplus \mathrm{E}_8(-1)^{\oplus 2} \oplus \langle -4 \rangle \cong \mathrm{Pic}(V)$ of rank-$19$. More generally, K3 surfaces polarized by the rank-19 lattice
	\begin{equation}
		\label{eq-M_N_lattice}
		\mathrm{M}_N := \mathrm{H} \oplus \mathrm{E}_8(-1)^{\oplus 2} \oplus \langle -2N \rangle
	\end{equation}
	for certain values of $N$ appear as mirrors to generic anticanonical K3 surfaces in Fano threefolds of Picard rank-1 \cite{MR2306141}. More details will be recalled in \S \ref{s-fano_mirror_symmetry}.
	\subsubsection{VHS for CY varieties}
	\label{sss-VHS}
	\par For any smooth $(n-1)$-dimensional CY variety $V$, the holomorphic $(n-1)$-form $\Omega$ determines a cohomology class $[\Omega]$ that generates $H^{n-1,0}(V)$. The Hodge decomposition on $H^{n-1}(V,\mathbb{C})$, polarized by intersection form, induces the usual Hodge filtration
	\begin{equation}
		\label{eq-hodge_filtration}
		\{0\} \subset F^{n-1} \subset F^{n-2} \subset \cdots \subset F^1 \subset F^0 = H^{n-1}(V,\mathbb{C})
	\end{equation}
	with $F^p = \oplus_{k \geq 0} H^{p+k,n-1-(p+k)}(V)$, so that $[\Omega]$ generates $F^{n-1}$ and that appropriate members of the filtration are orthogonal complements. The form $\Omega$ can be integrated over an integral basis $\{\Sigma_0,\dots,\Sigma_{k-1}\} \subset H_{n-1}(V,\mathbb{Z})/\mathrm{torsion}$, where $k=b_{n-1}$ is the middle Betti number of $V$, obtaining \emph{period integrals} $\omega_i$, $i=0,\dots,k-1$, given by
	\begin{equation}
		\label{eq-periods}
		\omega_i = \oint_{\Sigma_i} \Omega \, \equiv \langle \Sigma_i, \Omega \rangle \, ,
	\end{equation}
	where the latter notation denotes the Poincar\'e pairing. There is also a suitable interpretation when $V$ is an open CY variety. 
	\par There is an unobstructed complex structure deformation space $X_0$ of $\dim(X_0) = h^{n-2,1}(V)$ admitting the structure of a complex manifold, whose points $t \in X_0$ measure the choices inherent in (\ref{eq-hodge_filtration}); this is the Bogomolov-Tian-Todorov theorem \cite{MR915841,MR1027500}, see \cite{MR3592695} for results relevant to this article. Over $X_0$, there is a canonically flat, complex vector bundle $(\mathsf{H},\nabla) \to X_0$ whose fibres are the relevant complex cohomology groups $H^m(V_t,\mathbb{C})$. The flat connection $\nabla$ is called the Gauss-Manin connection, and $\mathsf{H}$ admits a local splitting according to the Hodge decomposition. Then the filtration (\ref{eq-hodge_filtration}) determines a sequence of holomorphic sub-bundles of $\mathsf{F}^p \subset \mathsf{H}$ and which obey Griffiths transversality with respect to $\nabla$, that is $\nabla \mathsf{F}^p \subseteq \mathsf{F}^{p-1}$. This determines the data of a polarized VHS. We assume that there is a compactification $X_0 \hookrightarrow X = X_0 \cup S$ obtained by including boundary components $S$ of $X_0$, stratified by components of dimension strictly less than that of $X_0$, that provide a measure of how $V$ degenerates in various limiting complex structures and acquires certain well-behaved singularities. Since $X_0 \subset X$ is rarely simply connected, there is a topological monodromy group $\Gamma$ that acts on the integral cohomology lattice $H^{n-1}(V_t,\mathbb{Z})/\mathrm{torsion} \subset H^{n-1}(V_t,\mathbb{C})$. There is a classifying space $\mathscr{D}$ of polarized Hodge structures of weight $(n-1)$ called the period domain, which can naturally be given the structure of a complex manifold. Accordingly, there is a well defined (holomorphic) map $X_0 \to \Gamma \backslash \mathscr{D}$ called the period map, which encodes the equivalent data of the Gauss-Manin connection. This in turn determines the VHS. See again Griffiths \cite{MR229641,MR233825,MR282990}. We assume as well that there is a limiting mixed Hodge structure defined over all of $X$ that extends the original VHS in a natural way.
	\par Inherently, the holomorphic $(n-1)$-form $\Omega = \Omega_t$ and the corresponding periods $\omega_i=\omega_i(t)$ depend on the complex structure $t \in X_0$. Suppose in addition that there is a polarized family of smooth, projective $(n-1)$-dimensional CY varieties $\pi : \mathcal{V}_0 \to X_0$ with general fibre $V_t := \pi^{-1}(t)$ whose complex structure is represented by $t \in X_0$. Moreover, assume that there is an extension $\mathcal{V}$ of $\mathcal{V}_0$ to all of $X$ where the varieties over the boundary strata $S$ acquire singularities; the Type III boundary components are isolated points where the variety degenerates to a union of hyperplanes and correspond to large complex structure limits in the context of mirror symmetry. The periods $\omega_i = \omega_i(t)$ become local analytic functions and determine a period sheaf $\Pi \to X_0$, whose rank is the number of linearly independent period integrals and carries the associated monodromy representation of $\Gamma$. The variation of the periods $\omega_i(t)$ governs crucial geometric data of the family $\mathcal{V}$ (the cycles $\Sigma_i$ can always be chosen to be locally constant in $t$), as the Poincar\'e pairing (\ref{eq-periods}) is compatible with the Gauss-Manin connection
	\begin{equation}
		\label{eq-period_derivative}
		\partial_t \, \omega_i(t) = \partial_t \, \langle \Sigma_i, \Omega_t \rangle \, = \langle \Sigma_i, \nabla_{\partial_t}\Omega_t \rangle \, ,
	\end{equation} 
	and underlies the full VHS on  $(\mathcal{H},\nabla)$. In the case $\mathrm{dim}(X)=1$, it follows from above  that there is a Fuchsian operator
	\begin{equation}
		\label{eq-PF_op}
		\mathcal{L}_n = P_0(\theta) + P_1(\theta) + \cdots + t^k P_k(\theta)
	\end{equation}
	of order $n$ called the Picard-Fuchs operator that annihilates the periods $\omega_i(t)$ of the family $V_t$, and underlies a VHS $\mathsf{H}$ of weight $n-1$ which we take to have nonzero Hodge numbers\footnote{For the $n-1=2$ case of K3 surfaces, we are assuming that $\mathrm{rank}(\mathrm{NS}(V))= 19 = 22-3$ so that $V$ has a 1-dimensional complex deformation space.} $h^{n-1,0} = h^{n-2,1} = \cdots = h^{0,n-1} = 1$. The operator $\mathcal{L}_n$ has regular singularities that coincide with the finite number of points in $S$, and measures the infinitesimal period relations of the family $\mathcal{V}_0$. Here $\theta = t \frac{d}{dt}$ is the logarithmic derivative and each $P_j(\theta)$ is a polynomial of order $n$. We will always assume that $X \cong \mathbb{P}^1$ and that $X_0$ is uniformized by the upper half plane $\mathbb{H} \subset \mathbb{C}$; thus $|S| \geq 3$.
	\par From the considerations of mirror symmetry, $\mathcal{L}_n$ has maximally unipotent monodromy (MUM) at $t=0$ \cite{MR1191426,MR1265317}, so $P_0(\theta)=\theta^n$.\footnote{Up to projective equivalence, MUM at a regular singular point is equivalent to all roots of the indicial equation being equal.} The Picard-Fuchs operator is realized as the regularization via Fourier-Laplace transform of the A-side small quantum D-module on the mirror CY $\check{V}$ at the corresponding large radius limit. In Manin \cite[Ch. II, \S 1]{MR1702284}, this is known as the second structure connection. There is a unique (up to scale) monodromy invariant holomorphic section $[\Omega]$ of $\mathsf{F}^{n-1}$ whose corresponding period integral $\omega$ is annihilated by $\mathcal{L}_n$. The variation of period integrals of the compact variety $V$ can be utilized with the mirror map to compute all genus zero invariants of the mirror variety $\check{V}$, which are encoded in a generating function known as the Yukawa coupling \cite[\S 2]{MR1191426}
	\begin{equation}
		\label{eq-yukawa}
		\mathcal{Y}_{n-1}(t) := \int_{V_t} \Omega_t \wedge \frac{d^{n-1}}{dt^{n-1}}\Omega_t = \langle \Omega_t , \nabla_{\partial_t}^{n-1} \Omega_t \rangle \, .
	\end{equation} 
	This is of particular interest in the case of $n-1=3,4$ for CY threefolds and fourfolds, in part due to the triviality of such invariants for the lower dimensional CY varieties of elliptic curves and K3 surfaces. That is to say, there are suitable normalizations in which $\mathcal{Y}_1$ and $\mathcal{Y}_2$ are constant. Of course, there are no nonconstant maps from $\mathbb{P}^1$ to an elliptic curve, and the statement for K3 surfaces can be traced back to the necessary quadratic period relations satisfied by K3 periods \cite{Nagura1995MirrorSymmetry}.
	\par This can be seen as follows. Suppose that $(\check{V},V)$ form a mirror pair of CY $(n-1)$-folds as described above. The mirror map can be constructed from $\mathcal{L}_n$ as follows. At the large complex structure limit $t=0$, there is a canonical Frobenius basis of solutions $\omega_0,\dots,\omega_{n-1}$ of the form
	\begin{align}
		\label{eq-MUM_basis}
		\omega_0(t) &= h_0(t) \, , \\
		\omega_1(t) &= \frac{1}{2\pi i}\left(h_0(t) \log(t)+h_1(t)\right) \, , \nonumber \\
		\vdots & \nonumber \\
		\omega_{n-1}(t) &= \frac{1}{(2\pi i)^{n-1}}\left(\sum_{k=0}^{n-1} \binom{n-1}{k} h_k(t)\log^{n-1-k}(t) \right) \,, \nonumber 
	\end{align}
	where $h_0,\dots, h_{n-1}$ are holomorphic functions satisfying $h_0(0)=1, h_1(0)=h_2(0)= \dots = h_{n-1}(0) = 0$. Then one sets 
	\begin{equation}
		\label{eq-mirror_map}
		\begin{aligned}
			\tau :& \; X_0 \to \mathbb{H} \, , \\
			\tau(t) =& \, \frac{\omega_1(t)}{\omega_0(t)} \, ,\\
		\end{aligned}
	\end{equation}
	which serves as the truncated period morphism. The map $\tau$ can be extended to a neighborhood $U \subseteq X$ of $t = 0$ by setting $\tau(0) = i \infty \in \mathbb{H}^* = \mathbb{H} \cup \mathbb{Q} \cup \{i\infty\}$. This extended map is invertible in a neighborhood $\check{U} \subseteq \mathbb{H}^*$ of $i \infty$, and the inverse $t(\tau)$ is known as the (closed) mirror map, as $i \infty \in \mathbb{H}^*$ corresponds to the large radius limit in the complexified symplectic moduli space of $\check{V}$. It is convenient then to define $q = \exp \circ \, 2 \pi i\tau : U \to \mathbb{C}$, so that
	\begin{equation}
		\label{eq-q_coordinate}
		q(t) = \exp(2 \pi i\tau(t)) = \exp\left(2 \pi i \frac{\omega_1(t)}{\omega_0(t)} \right) = t \exp\left(\frac{h_1(t)}{h_0(t)}\right) 
	\end{equation} and $q(0)=0$. In cases relevant to mirror symmetry, we have $q(t) \in t \mathbb{Z}[[t]]$. We will also refer to the local inverse $t = t(q)$ as the mirror map. 
	\par It is also convenient to write the Picard-Fuchs operator $\mathcal{L}_n$ in successive powers of $\theta$ as 
	\begin{equation}
		\label{eq-alt_pf}
		\mathcal{L}_n = \widetilde{P}_0(t)+\widetilde{P}_1(t)(\theta+1) + \cdots + \widetilde{P}_n(t)(\theta+1)^n \, .
	\end{equation}
	after the normalization $\omega(t) \mapsto t \omega(t)$. An elementary consequence of Griffiths transversality is that  the function $ \widetilde{P}_{n-1}(t)$ determines the Yukawa coupling $\mathcal{Y}_{n-1}(t)$ (\ref{eq-yukawa}) by integrating a first order equation of Wronskian type
	\begin{equation}
		\label{eq-yukawa_ODE}
		\theta \cdot \mathcal{Y}_{n-1}(t) = -\frac{2}{n}  \widetilde{P}_{n-1}(t) \mathcal{Y}_{n-1}(t) \, .
	\end{equation} 
	This is equivalent to the well-known relation $\mathcal{Y}(t)$ satisfies \cite{MR1191426}. An appropriate Lambert series $q$-expansion of the (normalized) pullback of $\mathcal{Y}(t)$ by the inverse of the mirror map
	\begin{equation}
		\label{eq-yukawa_expansion}
		\widetilde{\mathcal{Y}}_{n-1}(\tau) := \frac{1}{\omega_0(\tau)^2}\left(\frac{1}{2\pi i}\frac{dt}{d\tau}\right)^{n-1}\mathcal{Y}_{n-1}(t(\tau)) = 1 + \sum_{k=1}^{\infty} k^j N_k \frac{q^k}{1-q^k}
	\end{equation} 
	encodes the degree $k$ instantons as $N_k \in \mathbb{Z}$, up to an overall integer multiple of (\ref{eq-yukawa_expansion}). Here $j=3,2$ for $\dim(V) = n-1=3,4$ \cite{MR2282972}. For elliptic curves, $\mathcal{Y}(t)$ is up to a constant equal to the Wronskian of the Picard-Fuchs operator and so (\ref{eq-yukawa_expansion}) is easily seen to be constant \cite[Prop. 3.1.9]{zhou2014}, \cite[\S 5.1]{MR3822918}; for lattice polarized K3 surfaces of Picard rank $\rho = 19$ it is elementary that $\omega_0 \omega_2 = \omega_1^2$, and using Griffiths transversality this quadratic period relation forces (\ref{eq-yukawa_expansion}) to again be constant \cite[\S 4.2]{Nagura1995MirrorSymmetry}. Thus for $N_k = 0$ for all $k \in \mathbb{N}$ in both cases. 
	\subsection{Local mirror symmetry}
	\label{ss-local_mirror_symmetry}
	\par Local mirror symmetry specializes the techniques of mirror symmetry to study the geometry of a Fano subvariety $F \subset \check{V}$ \cite[\S 1]{MR1797015}. The large complex structure limit $t=0$ of the mirror $V$ corresponds to the limit in which $\mathrm{vol}(F) \to 0$ in $\check{V}$. This situation can be effectively modeled by the geometry of an open CY variety $\check{V}^\circ = \mathrm{Tot}(E)$, where $E \to F$ is a holomorphic vector bundle such that $K_{\check{V}^\circ} \cong \mathcal{O}_{\check{V}^\circ}$ is trivial. A common such example relevant to our considerations is the total space of the canonical line bundle $K_{\mathbb{P}^n} = \mathrm{Tot}(\mathcal{O}(-(n+1)))$ over projective space $\mathbb{P}^n$, as in \S \ref{s-intro}. This open CY variety of dimension $n+1$ is known as local $\mathbb{P}^n$. It is a toric variety whose mirror can be constructed from the procedure of Batyrev and Borisov \cite{MR1416334}; the compact part of the mirror is the mirror of a generic anticanonical divisor $\check{V} \in |-K_{\mathbb{P}^n}|$. A similar story holds for other Fano varieties, as before we will be exclusively concerned with $n=2,3$, so that the total space forms a local CY variety of dimension $n+1=3,4$, respectively.
	\par Coates and Givental \cite{MR2276766} showed that the A-side quantum $D$-modules of subvarieties $\check{V} \subset F$ determined by the vanishing of generic sections of a holomorphic vector bundle $E \to F$ can be determined by relation to the ambient quantum $D$-module of $F$ by quantum Lefschetz; in turn, quantum Serre duality relates GW invariants of $F$ twisted by the (equivariant) inverse Euler class of $E$ to the local invariants of $\check{V}^\circ = \mathrm{Tot}(E)$. Iritani, Mann, and Mignon \cite{MR3519131} showed that there is a non-equivariant limit in which quantum Serre duality can be expressed as a duality between the local and Euler-twisted quantum D-modules. This is directly recognizable when $E = K_F$ is the canonical line bundle. The mirror of a $F$ is a Laurent polynomial $\mathcal{W} \to \mathbb{P}^1_t - \{0,\infty\}$, thought of as a holomorphic map known as the superpotential. The generic fibres $V_t = \mathcal{W}^{-1}(t)$ are CY varieties that are believed due to physical considerations to be mirror to the generic anticanonical divisor $\check{V} \in |-K_F|$. The points $t=0,\infty$ correspond to singular CY varieties, and $t = \infty$ corresponds to the large complex structure limit of the corresponding family $\mathcal{V}_0 \to \mathbb{P}^1_t - \{0,\infty\}$ of CY varieties. There is more generally a finite set of singular points $\{0,\infty\} \subset S$ in which the fibre $V_t$ becomes singular, and over the Zariski open set $X_0 = \mathbb{P}^1 - S$ the family $\mathcal{V}_0 \to X_0$ is smooth. The totality of this data is known as a Landau-Ginzburg (LG) model. 
	\par In low dimensions, this relation is known to hold, see, e.g., Doran et al \cite{MR4913468} for a detailed case by case treatment of Fano threefolds and their anticanonical K3s. For each deformation class of Fano threefold, the superpotential $\mathcal{W}$ is explicitly constructed, and it is shown that generic fibres are lattice polarized K3 surfaces that obey the Dolgachev-Nikulin-Pinkham relation described in \S \ref{ss-classical_mirror_symmetry}. Let $\mathsf{H}$ be a weight $(n-1)$ VHS of the mirror family $\mathcal{V}_0$ as in \S \ref{sss-VHS} with associated Picard-Fuchs operator $\mathcal{L}_n$. Relevant to our considerations, it is well-known that the VHS of the mirror $V$, for the Fano surfaces and threefolds considered in this article, can be given a so-called modular parameterization \cite{MR1420220,MR2306141}. This means that there is a modular group $\Gamma \subset \mathrm{PSL}(2,\mathbb{R})$ and a Hauptmodul(= modular function of weight zero) $h : \Gamma \backslash \mathbb{H}^* \to X_0$ such that the pullback $h^*\mathsf{H}$ is isogenous to a symmetric power of the standard VHS on the relative $H^1$ of the elliptic modular surface over the modular curve $\Gamma \backslash \mathbb{H}^*$. In particular, the holomorphic period $\omega_0$ becomes a weight-$(n-1)$ modular form for $\Gamma$ when pulled back by the Hauptmodul = mirror map. Details related to this construction for modular elliptic surfaces will be recalled in Appendix \ref{s-modular_DEs}. For LG-models, such a modular parameterization is only possible for $n=2,3$ because of the presence of both MUM and conifold points in the discriminant locus $S$ rules out an $\mathrm{SL}(2)$ Mumford-Tate / geometric monodromy / differential Galois group more generally \cite[\S 10.6]{MR2851155}.\footnote{We thank the referee for explaining this point.}
	\par The Batyrev-Borisov procedure yields the operator $\check{\mathcal{D}}_{n+1} = \mathcal{L}_n \circ \, \theta$, which annihilates the periods of the mirror variety of $K_{\mathbb{P}^n}$. It is known as the local (small) quantum D-module \cite[\S 6]{MR3519131}. Accordingly, there is a constant period, which coincides with the noncompact direction in the open variety. This operator is also MUM at $t=0$, and the mirror map $\tau_{\mathrm{loc}}$ may be computed in exactly the same manner as before (\ref{eq-mirror_map}, \ref{eq-q_coordinate}). It is easy to see in this case that 
	\begin{equation}
		\label{eq-local_mirror_map_0}
		\tau_{\mathrm{loc}}(t) = \frac{1}{2\pi i}\int_{t_0}^t \frac{\omega_0(s)}{s} \, ds 
	\end{equation}
	for a fixed choice $t_0 \in X_0$.
	\par In the context of this article, we will be interested primarily in the closely related operator
	\begin{equation}
		\label{eq-quantum_operator}
		\mathcal{D}_{n+1} := \theta \circ \mathcal{L}_n
	\end{equation}
	which we deem the \emph{quantum operator}, see \S \ref{s-fano_mirror_symmetry}. This extended operator is also MUM at $t = 0$, provided that $\mathcal{L}_n$ is suitably normalized \cite[\S 9]{MR4503993}, and moreover is closely related to the local operator $\check{\mathcal{D}}_{n+1}$. Assume that $\mathsf{H}$ admits a modular parameterization. It follows from Doran and Kerr \cite{MR2851155} that $\mathcal{D}_{n+1}$ underlies a mixed VHS (MVHS) $\mathsf{V}$ of the form 
	\begin{equation}
		\label{eq-MVHS}
		0 \to \mathbb{Q}(1) \to \mathsf{V} \to \mathsf{H} \to 0
	\end{equation}
	with dual sequence (after a Tate twist)
	\begin{equation}
		\label{eq-dual_sequence}
		0 \to \mathsf{H} \to \mathsf{V}^{\vee}(-n+1) \to \mathbb{Q}(-n) \to 0 \, .
	\end{equation}
	One choice of Picard-Fuchs operator for (\ref{eq-dual_sequence}) is the local operator $\check{\mathcal{D}}_{n+1}$. The latter extension arises from a higher cycle $\xi \in H^n(\mathcal{V}_0,\mathbb{Q}(n))$, and the associated extension regulator class $R \in \Gamma(X_0,\mathsf{H} \otimes \mathcal{O}_{X_0}/\mathsf{H} \otimes \mathbb{Q}(n))$ is annihilated by the local operator, $\check{\mathcal{D}}_{n+1}R=0$. Dually, the multivalued function $\nu := \langle R,\Omega \rangle$ is annihilated by the quantum operator, $\mathcal{D}_{n+1}\nu = 0$ precisely when $\mathcal{L}_n$ is normalized as in (\ref{eq-alt_pf}) for the computation of the Yukawa coupling. An elementary consequence of variation of parameters is that the higher normal function $\nu(t)$ is independent from the periods $\omega_0(t),\dots,\omega_{n-1}(t)$ annihilated by $\mathcal{L}_n$ (and hence by $\mathcal{D}_{n+1}$), and we may write $\nu(t) = \omega_n(t) \sim \log^n(t)$ to complete the Frobenius basis (\ref{eq-MUM_basis}) for the quantum operator $\mathcal{D}_{n+1}$.
	\subsection{Modularity and GW invariants of local threefolds}
	\label{ss-local_threefolds}
	\par Several properties of the higher normal function $\omega_{n}(t)$ and extension regulator class $R$ relevant to modularity and local mirror symmetry have been proven by Bloch, Kerr, and Vanhove \cite{MR3780269} and Doran and Kerr \cite{MR2851155}. Continuing in the spirit from the previous section (again recall $n=2,3$), from the VHS $\mathsf{H}$, choose a complementary locally $\nabla$-constant $\mathbb{Q}$-Betti basis $\Sigma_0,\dots,\Sigma_{n-1} \in H_{n-1}(V_t,\mathbb{Q})$ near the large complex structure limit $t=0$ for the family $\mathcal{V}_0$. Then setting $R_i := \langle \Sigma_i,R \rangle$, the local (exponentiated) mirror map is determined by 
	\begin{equation}
		\label{eq-local_mirror_map}
		Q(t) := \exp(2\pi i \tau_{\mathrm{loc}}(t)) \, ,
	\end{equation}
	where the local mirror map (\ref{eq-local_mirror_map_0}) is identified $\tau_{\mathrm{loc}}(t) \equiv \frac{1}{(2\pi i)^2}R_0(t)$. From the exponentiated mirror map $q(t)$ in (\ref{eq-q_coordinate}), the parameter $t$ can be eliminated, and then we have 
	\begin{equation}
		\label{eq-local_mirror_map_derivative}
		q \frac{d}{dq} \log Q = (2\pi i) \frac{d}{d\tau} \tau_{\mathrm{loc}} \, .
	\end{equation}
	Moreover, for LG models admitting a modular parameterization for a modular group $\Gamma_m$ (here $m \in \mathbb{N}$ labels the group as e.g., $\Gamma_m = \Gamma_0(m)$, see Appendix \ref{s-modular_DEs}), this quantity is closely related to the normalized higher normal function 
	\begin{equation}
		\label{eq-eisenstein_potential}
		\mathscr{E}_{n,m} := \langle R, \hat{\Omega} \rangle = \frac{\omega_n}{\omega_0}
	\end{equation}
	where $\hat{\Omega} = \frac{\Omega}{\omega_0}$. We refer to $\mathscr{E}_{n,m}$ as a \emph{rank-$n$ Eisenstein potential} for $\Gamma_m$. We have the following result due to Doran and Kerr which explains our naming convention and plays a crucial role in this article. 
	\begin{theorem}[\cite{MR2851155}, Prop. 9.1]
		\label{thm-eisenstein}
		The rank-$n$ Eisenstein potential $\mathscr{E}_{n,m}$ is an Eichler integral of a $\mathbb{Q}$-Eisenstein series of weight-$(n+1)$ as
		\begin{equation}
			\label{eq-eichler_integral}
			(2\pi i)^n\frac{d^{n}}{d\tau^{n}}\mathscr{E}_{n,m} = (2\pi i)\frac{d}{d\tau} \tau_{\mathrm{loc}} = E_{n+1,m} \, .
		\end{equation}
	\end{theorem}
	\par Here, it is meant (for $n=2,3$) that $E_{n+1,m}$ defines a weight-$(n+1)$ modular form of Eisenstein type for $\Gamma_m$ and has constant terms that are rational at the cusps. These modular forms appear naturally in mirror symmetry for Fano threefolds of rank-1. It is crucial for us that such modular forms are Eisenstein series. For reasons that are explained below, it is quite natural that expressions like those in (\ref{eq-eichler_integral}) are modular forms of the appropriate weight, as can be seen below by their natural factorization (\ref{eq-eisenstein_factorization}). The fundamental reason that Eisenstein series appear is that these forms are realized as the pullback of a motivic cohomology class under a map from a Kuga variety to an LG-model, and the cycle on the Kuga variety is necessarily an Eisenstein symbol. This is explained by Doran and Kerr \cite[\S 10]{MR2851155}.
	\par Each deformation class of such a Fano threefold is determined in part by an integer invariant $N \in \{1,2,3,4,5,6,7,8,9,11\}$ called the level, see \S \ref{ss-Fano_MS_&_QDE}. Golyshev and Zagier showed in \cite[Equation (1.14)]{MR3462676} that for each $N > 1$, $E_{4,N}(\tau)$ is of weight-4 for $\Gamma_0(N)^+$ and can be expressed as a rational linear combination of Eisenstein series as
	\begin{equation}
		\label{eq-eisenstein_expansion}
		E_{4,N}(\tau) = A_N \sum_{M | N} M^2 h_M G_4(M\tau) \, ,
	\end{equation}
	where $G_4(\tau)$ is the Hecke-normalized Eisenstein series of weight-4, i.e., 
	\begin{equation*}
		G_4(\tau) = \frac{1}{240}E_4(\tau) = \frac{1}{240} + \sum_{n=1}^{\infty} \frac{n^3 q^n}{1-q^n} \, .
	\end{equation*}
	The usual Eisenstein series $E_4(\tau)$ is given in Equation (\ref{eq-eisenstein_series_Ek}). Similar expressions appear in Doran and Kerr \cite[\S 10.5]{MR2851155}. Here $A_N \in \mathbb{Q}$ is an overall scale factor, and the coefficients $h_M \in \mathbb{Q}$ enjoy the antisymmetry property $h_M = - h_{N/M}$. Specific expressions will be provided in \S \ref{s-fano_mirror_symmetry}. In fact, these $\mathbb{Q}$-Eisenstein series satisfy a natural factorization property: a straightforward computation using the modular parameterization of the Picard-Fuchs operator shows that $E_{n+1,m}$ is (up to an overall rational multiple) the product of (meromorphic)
	modular forms
	\begin{equation}
		\label{eq-eisenstein_factorization}
		E_{n+1,m}(\tau) = \omega_0(\tau) \frac{t^{\prime}(\tau)}{t(\tau)}  \, ,
	\end{equation}
	where $t^{\prime} = \frac{1}{2\pi i}\frac{dt}{d\tau}$. See also \cite[\S 10.2]{MR2851155} or Theorem \ref{thm-modular_ODEs} for more details.
	\par In the case of $n=2$, this construction plays an important role in mirror symmetry of local threefolds as follows. Let $\mathcal{X}$ be the LG model of a toric Fano surface $F$. It was conjectured by Hosono \cite{MR2282969} and proven by Bloch, Kerr, and Vanhove \cite[Thm. 6.1]{MR3780269} that $R_1$ is the local genus zero GW potential, that is
	\begin{equation}
		\label{eq-R_1_GW}
		R_1(\tau_{\mathrm{loc}}) = p(\tau_{\mathrm{loc}}) - \sum_{k=1}^{\infty} k n_k Q^k 
	\end{equation}
	for certain genus zero GW invariants $n_k$ of the local threefold $K_F$, and $p$ a quadratic polynomial that is obtained after restricting to the anti-canonical line in the complexified symplectic moduli. A consequence of this is the following, due to Doran and Kerr \cite[\S 5]{MR2851155}. Differentiating twice with respect to $\tau_{\mathrm{loc}}$, one obtains 
	\begin{equation}
		\label{eq-local_instantons_n=2}
		(2\pi i) \frac{\mathcal{Y}_1(t(Q))}{\omega_0(t(Q))^3} = (2\pi i)^2 \frac{d^2 R_1}{d R_0^2} = r - \sum_{k=1}^{\infty} k^3 n_k Q^k
	\end{equation}
	where $r \in \mathbb{N}$ and as previously $\mathcal{Y} = \langle \Omega,\theta \Omega \rangle$. Up to a change in normalization, one obtains an instanton expansion like the one in (\ref{eq-yukawa_expansion}) for $j=3$, computing the degree $k$ genus zero GV invariants $N_k$ for $K_F$. This is in fact quite natural from the expression for $\widetilde{\mathcal{Y}}_{1}$ in (\ref{eq-yukawa_expansion}). Notice further that this expression is a natural inversion of (\ref{eq-local_mirror_map_derivative}) as follows: we have
	\begin{equation}
		\label{eq-inversion_eq1}
		\frac{d R_1}{d R_0} = \frac{\theta R_1}{\theta R_0} = \frac{\omega_1}{\omega_0} = \tau \, ,
	\end{equation}
	and thus
	\begin{equation}
		\label{eq-inversion_eq2}
		(2\pi i)^2 \frac{d^2 R_1}{d R_0^2} = (2\pi i) \frac{d}{d\tau_{\mathrm{loc}}}\tau =  Q\frac{d}{dQ} \log q \, .
	\end{equation}
	Our goal for the remainder of this article is to provide evidence that the analogous expansions from (\ref{eq-inversion_eq2}) for $n=3$, $j=2$ provide instanton computations of local CY fourfolds $K_F$ over rank-1 Fano threefolds. This discussion occurs in \S \ref{s-genus_zero_invariants}.
	%
	%\par There is a classical geometric interpretation of second order Fuchsian equations with three or four singular points on $\mathbb{P}^1$ in terms of pencils $\pi : \mathcal{E} \to \mathbb{P}^1$ of elliptic curves is known as Gauss-Schwarz theory \cite{MR2514149}, which we describe now. 
	%	
	
	\section{Modularity of Fano mirrors}
	\label{s-fano_mirror_symmetry}
	\par It is well established in the literature that mirror symmetry of the Fano surfaces under consideration in this article is governed by the geometry of modular elliptic surfaces, see, e.g., \cite{MR1486340,MR1468700,MR1732409} as well as many of the sources in the previous section, such as \cite[\S 10.4]{MR2851155}. These are the canonical elliptic surfaces $\pi : \mathcal{E} \to X$ over a modular curve $X = \Gamma \backslash \mathbb{H} \cup \{\mathrm{cusps}\}$, where $\Gamma \subseteq \mathrm{PSL}(2,\mathbb{Z})$ is a congruence subgroup, and constitute important examples of extremal elliptic surfaces \cite{MR0867347}. Then the Picard-Fuchs operator $\mathcal{L}_2$ annihilating the periods of the smooth elliptic fibres admit a modular parameterization as described in \S \ref{ss-local_mirror_symmetry}; details are provided in Appendix \ref{s-modular_DEs} and Appendix \ref{ss-wm}. Similarly, the modularity of the mirrors rank-1 Fano threefolds has been known since the work of Dolgachev \cite{MR1420220} and in particular Golyshev \cite{MR2306141}. Roughly, this can be traced geometrically to the mirror K3s admitting a canonical $\mathrm{M}_N$-polarization, whose coarse moduli space is the modular curve $X_0(N)^+$, together with the necessary quadratic period relations. Such K3s are well-known to be Shioda-Inose partners of the Kummer surface of the product of two elliptic curves that are $N$-isogenous.
	\subsection{Mirror symmetry of Fano varieties and quantum differential equations}
	\label{ss-Fano_MS_&_QDE}
	
	\par  In the context of mirror symmetry of Picard rank-1 Fano threefolds, Golyshev \cite{MR3415413} introduced a class D$n$ of $n^{th}$-order Fuchsian ordinary linear differential operators $\mathcal{L}_n$. There, the differential operators are relevant to the (small) quantum cohomology of smooth Fano varieties. The operators were then studied in depth by Golyshev and Stienstra \cite{MR2346574}. A lengthy, precise definition of such operators is not needed here. However, let us summarize some relevant aspects: each determinantal D$n$ operator $\mathcal{L}_n$ is defined in a combinatorial way via determinants from an $(n + 1) \times (n +1)$-matrix $A = [a_{i j}]$ in such a way that $\mathcal{L}_n$ always arises as the rightmost factor of a Fuchsian operator $\mathcal{D}_{n+1}$ of order $(n+1)$, given by
	\begin{equation}
		\label{eq-determinantal_eq}
		\mathcal{D}_{n+1} = \theta \circ \mathcal{L}_n \, \, .
	\end{equation}
	See Golyshev \cite[Definition 2.10, \S 5.6]{MR2306141} for the precise definition and relevant matrices $A=[a_{i j}]$. 
	\par Determinantal D3 operators were studied by Golyshev in the context of mirror symmetry of the precisely seventeen algebraic deformation classes of rank-1 Fano threefolds, i.e., algebraic threefolds with ample anticanonical class such that the Picard lattice has rank-1 that were classified by Iskovskikh \cite{MR0463151,MR0503430}. We will not need a geometric description of these varieties and defer the interested reader to Iskovskikh's work; see also Golyshev \cite[\S 6.3]{MR2306141} or Coates et al \cite{MR3470714}. For each deformation class $F$, the corresponding operator $\mathcal{L}_3$ was shown to admit a modular parameterization as in \S \ref{ss-local_mirror_symmetry} for certain congruence subgroups $\Gamma \subseteq \mathrm{SL}(2,\mathbb{R})$, and to be the Picard-Fuchs operator of the pencil $\pi \colon \mathcal{V}_0 \to X_0 \subset \mathbb{P}^1$ of K3 surfaces of Picard rank-19 that are  Dolgachev-Nikulin-Pinkham mirror \cite{MR1420220} to the anticanonical linear system on $F$. Moreover, these operators were shown by Golyshev \cite[\S 2.10]{MR2306141} to have precisely the normalization described in (\ref{eq-alt_pf}). It follows then that the operator $\mathcal{D}_{n+1}$ is the associated quantum operator from (\ref{eq-quantum_operator}). Modular D3 operators and the associated quantum operators played a crucial role for Golyshev and Zagier \cite{MR3462676} in their proof of the Gamma conjecture for the associated rank-1 Fano threefolds.
	\par Moreover, modular D2 operators $\mathcal{L}_2$ were shown by Golyshev and Vlasenko to agree precisely with the modular second-order operators satisfying Ap\'ery like recurrences like those that arise in irrationality proofs for $\zeta(2)$, discovered by Zagier in \cite{MR2500571}. These are well-known to be precisely the Picard-Fuchs operators of the canonical elliptic fibrations over rational modular curves for congruence subgroups of $\mathrm{PSL}(2,\mathbb{Z})$, whose fibres also arise as mirrors of anticanonical elliptic curves in del Pezzo surfaces. For example, the canonical elliptic fibration $\pi : \mathcal{E} \to X_0(3) \cong \mathbb{P}^1$ serves as the mirror space for $\mathbb{P}^2$, in the sense that $t = 0 \in X_0(3)$ as a large complex structure limit for the family of anticanonical cubic curves in $\mathbb{P}^2$ \cite{MR1732409,MR4301054,MR3822918}. Such elliptic fibrations can be written as a Weierstrass model 
	\begin{equation}
		\label{eq-RES_wm}
		\mathcal{E}_t := \pi^{-1}(t) \; : \quad	y^2 = 4x^3 - g_2(t)x - g_3(t)
	\end{equation}
	where $t$ is a Hauptmodul of the modular curve and $g_2(t),g_3(t)$ are certain polynomials of degree $\leq 4,6$ respectively. Defining expressions are well-known throughout the literature, and we collect them in Table \ref{table-modular_elliptic_data}. These are extremal rational elliptic surfaces that have only three or four singular fibres, and were classified by Miranda and Persson \cite{MR0867347}. For us, the relevant modular curves are $X_0(m)$, $m=3,4,6,8,9$, and $X_1(5)$. Moreover, the modular elliptic surfaces over $X_0(8), \, X_0(9)$ can be realized as pullbacks of those over $X_0(4)$, $X_0(3)$ by the map $t \mapsto t^k$ with $k=2,3$ respectively, so it is sufficient to consider only those listed in Table \ref{table-modular_elliptic_data}.
	The holomorphic 1-form $\Omega_t = dx/y$ can be integrated over a basis of 1-cycles on each smooth fibre as in (\ref{eq-periods}), and near the large complex structure limit $t=0$, there is a canonical holomorphic period $\omega_0(t)$ that becomes a weight-1 modular form for the corresponding congruence subgroup. The weight-3 modular form $E_{3,m}$ of Eisenstein type is determined from the modular data as in (\ref{eq-eisenstein_factorization}). We list modular data of each family in Table \ref{table-modular_elliptic_data}, as well as the associated Picard-Fuchs operators and quantum operators in Table~\ref{table-modular_pf_elliptic}. The Picard-Fuchs operator can be computed from the data of the Weierstrass model as is recalled in \ref{eq-elliptic_pf}.
	\par Let us provide some more details for Fano threefolds: a deformation class $F=F_{N,d}$ of a rank-1 Fano threefold is determined by two numerical invariants, the index $d$ and the level $N$. Each class in $H^2(F,\mathbb{Z})$ is algebraic, so we have $H^2(F,\mathbb{Z}) \cong \mathbb{Z}$. The index is given by $d = \left[H^2(F,\mathbb{Z}) \, : \, \mathbb{Z}c_1\right]$, where $c_1$ is the anticanonical class, and the level is $N = \frac{1}{2d^2} \left\langle c_1^3,[F] \right\rangle$. Both are easily shown to be positive integers. The seventeen possible combinations of invariants $(N, d)$ determined by Iskovskikh are as follows:
	\begin{align}
		\label{eq-degree_level_invariants}
		(1,1), \; \; (2,1), \; \; (3,1), \; \; (4,1), \; \; & (5,1), \; \; (6,1), \; \;  (7,1), \; \; (8,1), \; \; (9,1), \; \; (11,1), \nonumber \\
		(1,2), \; \; (2,2), \; \; (3,2), \; \; & (4,2), \; \; (5,2), \; \; (3,3), \; \; (2,4)\, . 		
	\end{align}
	Each corresponding Picard-Fuchs operator $\mathcal{L}_3 = \mathcal{L}_{3,N}$ is normalized to have a MUM point at $t=0$ . In fact, each $(N, d)$-modular family is the pullback of the twisted symmetric square of the universal elliptic curve over the genus-zero modular curve $X_0(N)^+$ by a cyclic covering of degree $d$. Certain modifiers must be added to this sentence for the $N=1$ case, which we will not concern ourselves with. It then follows from Theorem \ref{thm-modular_ODEs} that $\mathcal{L}_{3,N}$ is projectively equivalent to the symmetric square $\mathcal{L}_2^{\otimes 2}$, where $\mathcal{L}_2$ is the second order Picard-Fuchs operator of an elliptic pencil, whose monodromy group is realized by the symmetric square representation of $\Gamma_0(N)^+$ in $\mathrm{SO}(2,1)$. It follows from Dolgachev \cite[Thm. 7.1]{MR1420220} that the mirror K3 pencils are closely connected with K3 surfaces polarized by the rank-19 lattice $\mathrm{M}_N$ from (\ref{eq-M_N_lattice}). In this section we focus on pairs $(N, d)$ in (\ref{eq-degree_level_invariants}) with $N \in \{2,3,4,5,6\}$ and $d=1$ as they can naturally be linked with the modular elliptic surfaces in Table \ref{table-modular_elliptic_data} via the twist description we will describe shortly; analogous results can be obtained for the remaining pairs by rational pullback. A specific instance of this will arise in \S \ref{s-genus_zero_invariants}. The relevant third order Picard-Fuchs operators are given in Table~\ref{table-modular_pf}; they were determined by Golyshev in \cite[\S 5.8]{MR2306141}. 
	\par The existence of a canonical $\mathrm{M}_N$-polarization will allow us to provide Weierstrass models for K3 elliptic fibrations, and provides an intrinsic description of these mirror K3 pencils using the twist construction of Doran and Malmendier  \cite{MR4069107} applied to the modular elliptic surfaces in Table \ref{table-modular_elliptic_data}; see Appendix \ref{ss-twisted_wm}. In short, the twist construction is determined by the generalized functional invariant, which is encoded as a triple $(i,j,\alpha)$ for $i,j \in \mathbb{N}$ and $\alpha \in \{\frac{1}{2},1\}$ labeling a 1-parameter family of degree $i+j$ covering maps $\mathbb{P}^1 \to \mathbb{P}^1$.  For homogeneous coordinates $[v_0 : v_1] \in \mathbb{P}^1$, the relevant family of maps parameterized by $\mathsf{t} \in \mathbb{P}^1 - \{0,1,\infty\}$ is given by 
	\begin{equation}
		\label{eqn-generalized_functional_invariant0}
		[v_0,v_1] \mapsto [c_{ij}v_1^{i+j}\mathsf{t} : v_0^i(v_0+v_1)^{j}],
	\end{equation}
	for a certain constant $c_{i j} \in \mathbb{C}^\times$. The significance of the generalized functional invariant is that it allows for the construction of towers of iterated CY fibration structures on a CY $n$-fold can be obtained starting from a given elliptically fibred CY $(n-1)$-fold pencil. Moreover, they showed how the new holomorphic period and Picard-Fuchs operator of the resulting family can be expressed in terms of Hadamard products of lower rank hypergeometric functions and operators, see Equation~\ref{eq-Hadamard_product}.
	\par It was shown in \cite[Lem. 6.9, \S 10.2]{MR4069107} that elliptic fibrations on $\mathrm{M}_N$-polarized K3 pencils can be obtained by the generalized functional invariant $(i,j,\alpha)=(1,1,1)$ from the canonical elliptic fibrations over the modular curves in Table \ref{table-modular_elliptic_data} for $m=N=3,4,5,6$. The twisted model has the form
	\begin{equation}
		\label{eq-(1,1,1)-twisted_wm}
		\scalemath{0.9}{y^2=4 x^3-g_2\left(-\frac{t}{ u(u+1)}\right)(u(u+1))^4 x-g_3\left(-\frac{t}{ u(u+1)}\right)(u(u+1))^6,}
	\end{equation}
	where $g_2,\,g_3$ are the original coefficients in the Weierstrass model (\ref{eq-RES_wm}) found in Table ~\ref{table-modular_elliptic_data}. Here we have normalized the twist variable to be $\mathsf{t}=4t$ from Equation (\ref{eqn-generalized_functional_invariant0}) to ensure integrality of the holomorphic K3 period. In this case, the relevant hypergeometric function that appears as the twist function in the Hadamard product is
	\begin{equation}
		\label{eq-1F0_0}
		_1F_0 \left(\frac{1}{2} \, \Bigg{|} \, 4t\right) = \frac{1}{\sqrt{1-4t}} = \sum_{k=0}^{\infty} \binom{2k}{k} \, t^k \, 
	\end{equation}
	This holomorphic function is annihilated by the differential operator 
	\begin{equation}
		\label{eq-1F0_op0}
		\mathcal{L}_1 = t(2+4\theta)-\theta \, ,
	\end{equation}
	and it follows that the Picard-Fuchs operator $\mathcal{L}_{3,N} = \mathcal{L}_1 \star \mathcal{L}_{2,N}$ of the $\mathrm{M}_N$-polarized K3 decomposes as the Hadamard product of $\mathcal{L}_1$ with the corresponding elliptic Picard-Fuchs operator $\mathcal{L}_{2,N}$ from Table \ref{table-modular_pf_elliptic}. 
	\par The $\mathrm{M}_2$-polarized family, mirror to the family quartic K3s in $\mathbb{P}^3$, can be obtained in a similar way from the modular curve $X_0(2)$ and generalized functional invariant $(i,j,\alpha)=(1,1,1)$, but we choose as in \cite[Lem. 7.3]{MR4069107} to use the generalized functional invariant $(i,j,\alpha) = (3,1,1)$ and pullback the modular elliptic surface over $X_0(3)$. This links the mirror of $\mathbb{P}^2$ with the mirror of $\mathbb{P}^3$ via the twist construction and iterative CY fibration structures. We obtain the minimal Weierstrass model
	\begin{equation}
		\label{eq-M_2_wm}
		\scalemath{.9 }{y^2=4 x^3-g_2\left(-\frac{t}{ u^3(u+1)}\right)(u(u+1))^4 x-g_3\left(-\frac{t}{ u^3(u+1)}\right)(u(u+1))^6,}
	\end{equation}
	where $g_2, \, g_3$ are the corresponding coefficients in Table \ref{table-modular_elliptic_data}. We have normalized the twist variable to be $\mathsf{t}=\frac{4^4}{3^3}t$ from Equation (\ref{eqn-twisted_weierstrass}); this ensures integrality of the holomorphic K3 period. It follows that Equation \ref{eq-M_2_wm} defines a universal family of $\mathrm{M}_2$-polarized K3 surfaces over a $1$-dimensional coarse moduli space constructed by Dolgachev \cite[Thm. 7.1]{MR1420220}; here the (compactified) moduli space is the rational modular curve $X_0(2)^+ \cong \mathbb{P}^1$, and the family (\ref{eq-M_2_wm}) is parameterized by a Hauptmodul $t$ given in Table \ref{table-K3_modular}. The corresponding Picard-Fuchs operator $\mathcal{L}_{3,2} = \mathcal{L}_3 \star \mathcal{L}_{2,3}$ decomposes under the Hadamard product with the third order hypergeometric operator
	\begin{equation}
		\label{eq-3F2_twist_op}
		\mathcal{L}_3 =  3\theta(3\theta-1)(3\theta-2)-8t(2\theta+1)(4\theta+1)(4\theta+3)
	\end{equation}
	that annihilates the hypergeometric function
	\begin{equation}
		\label{eq-3F2_twist_fn}
		_3F_2 \left(\begin{array}{c}
			\frac{1}{4} \, , \, \frac{2}{4} \, , \, \frac{3}{4} \\
			\frac{1}{3} \, , \, \frac{2}{3}
		\end{array} \Bigg {|} \, \frac{4^4}{3^3}t \right) = \sum_{k=0}^{\infty} \binom{4k}{k} t^k \, ,
	\end{equation}
	see Appendix \ref{eq-Hadamard_product}. 
	\par A particularly interesting such example is provided by the so-called Ap\'ery pencil $Y_t$ of K3 surfaces, represented in affine coordinates by the equation
	\begin{equation}
		\label{eq-apery_family}
		Y_t \colon \quad	1-(1-x y) z-t x y z(1-x)(1-y)(1-z)=0 \, .
	\end{equation}
	This family of K3 surfaces is generically of Picard rank $\rho=19$, and was utilized by Beukers and Peters \cite{MR0749676} in their ``geometrization'' of the irrationality proof of $\zeta(3)$. It follows from Golyshev \cite[\S 6.3]{MR2306141} that this family is mirror dual to the family of anticanonical K3 surfaces in the Fano threefold $V_{12}$, which is realized geometrically in the Iskovskikh classification as a section of the orthogonal Grassmannian $O(5,10)$ by a codimension-7 plane. Accordingly, the third order Picard-Fuchs operator is the operator $\mathcal{L}_{3,6}$ given in Table \ref{table-modular_pf}. 
	\par The Ap\'ery family does not admit an $\mathrm{M}_6$-polarization, but is closely related to a family that does. The family in question is represented by the affine equation
	\begin{equation}
		\label{eq-fermi_surface}
		\mathscr{Y}_s \; : \quad	X+\frac{1}{X}+Y+\frac{1}{Y}+Z+\frac{1}{Z}=s+\frac{1}{s} \;,
	\end{equation}
	with $s \in \mathbb{P}^1 - \{0,\infty\}$, and is called the Fermi pencil. It was shown in \cite[Thm. 1]{MR1034260} that $\mathscr{Y}_s$ represents a pencil of K3 surfaces of Picard rank-19 that generically admit a canonical $\mathrm{M}_6$-polarization. By \cite[Rmk. 4]{MR1034260}, the Beukers-Peters K3 surface $Y_t$ in Equation (\ref{eq-apery_family}) is obtained by the resolution of the quotient of the Fermi pencil $\mathscr{Y}_s$ in (\ref{eq-fermi_surface}) by the involution $(X,Y,Z,s) \mapsto (-X,-Y,-Z,-s)$. Hence, it follows that the periods of $Y_t$ lift to periods of (\ref{eq-fermi_surface}), up to a rational change of coordinates in loc. cit.. Alternatively, this relation can be seen at the level of special functions by applying a sequence of identities for the Heun function \cite[Rmk. 10.10]{MR4069107}. An elliptic fibration on (\ref{eq-fermi_surface}) realizing the $\mathrm{M}_6$-polarization can be obtained by applying the twisted model of (\ref{eq-(1,1,1)-twisted_wm}) to the Weierstrass data in Table \ref{table-modular_elliptic_data} for $X_0(6)$. It follows from the discussion above (see also \cite[\S 8, Ex. $\mathscr{C}$]{MR0783555} and \cite[\S 3]{MR1034260}) that the Picard-Fuchs operator for the twisted model decomposes as a Hadamard product $\widetilde{\mathcal{L}}_{3,6} = \mathcal{L}_1 \star \mathcal{L}_{2,6}$, and is equivalent to the Picard-Fuchs operator on the Fermi pencil. This operator is found under the $N = \widetilde{6}$ entry in Table \ref{table-modular_pf}. It is also of interest to note that the Fermi pencil in (\ref{eq-fermi_surface}) and its 2:1 cover of the Ap\'ery family (\ref{eq-apery_family}) was utilized by Kerr in \cite{kerr_motivic_2020} to produce $lcm^3$ bounds on the denominators of Ap\'ery's second series in the irrationality of $\zeta(3)$.
	\par Golyshev and Vlasenko observed in \cite[\S 8]{MR3415413} that up to rational pullback, the modular D3 operators for all families we have associated with generalized functional $(i, j,\alpha)=(1,1,1)$ could be obtained from modular D2 operators by convolution, in perfect agreement with our computations. Thus, we have geometrized their results, and extended in one case to other generalized functional invariants. In Table \ref{table-K3_modular}, we provide modular data of the relevant $\mathrm{M}_N$-polarized K3 families, as well as for the Ap\'ery family. In Table \ref{table-modular_pf} we list the Picard-Fuchs operator, twist data, and quantum operators. Similar to the Ap\'ery operator $\mathcal{L}_{3,6}$, the twisted operators in the $N = 8,9$ cases are distinct from the corresponding operators found in \cite[\S 5.8]{MR2306141}, but are linked to those by similar special function identities in \cite[Rmk. 10.10]{MR4069107}. We omit these computations as they are very similar in spirit to what we have done in this section. The $N=7,11$ families are not readily obtained via the twist construction; the modular data and Picard-Fuchs operators can be found in \cite{MR2306141}, and weight-4 $\mathbb{Q}$-Eisenstein series computed analogously as we do here.
	\section{Genus zero invariants of local \texorpdfstring{$\mathbb{P}^2$}{P^2} and local \texorpdfstring{$\mathbb{P}^3$}{P^3}}
	\label{s-genus_zero_invariants}
	\par In this section we demonstrate on an experimental basis that the discussion in \S \ref{ss-local_threefolds} for computing genus zero GV invariants for local CY threefolds can be generalized to compute such invariants for local CY fourfolds. First, we define the local CY fourfold analogue of the data given in \S \ref{ss-local_threefolds}.
	\subsection{Functional relations between local and closed mirror maps}
	\label{ss-functional_relations}
	Recall from that section (\ref{eq-local_instantons_n=2}, \ref{eq-inversion_eq2}) that the expression 
	\begin{equation}
		\label{eq-local_threefold_Q_expansion0}
		Q \frac{d}{dQ} \log q = r - \sum_{k=1}^{\infty} k^3 n_k Q^k
	\end{equation}
	encodes the degree $k$ genus zero GW invariants $n_k$ of an appropriate local threefold. This essentially proves an observation of Stienstra \cite{MR2282958}, that the associated Lambert series expansion
	\begin{equation}
		\label{eq-local_threefold_Q_expansion1}
		Q \frac{d}{dQ} \log q = r + \sum_{k=1}^{\infty} k^3 \widetilde{\mathsf{N}}_k \frac{Q^k}{1-Q^k}
	\end{equation}
	recovers the instanton expansion of said space, that the $\widetilde{N}_k$ are the degree $k$ genus zero GV invariants.\footnote{Note the change in normalization. This also differs by a factor of $-1$ from Stienstra \cite{MR2282958}.} Let us label the invariants as $\widetilde{\mathsf{N}}_{k}=\widetilde{\mathsf{N}}_{m,k}$.
	\par This expression (\ref{eq-local_threefold_Q_expansion1}) can be understood as an inversion of the weight-3 $\mathbb{Q}$-Eisenstein series $E_{3,m}(\tau)$ that are given in Table \ref{table-modular_elliptic_data}. Let us write such a series as a $q$-Lambert series itself,
	\begin{equation}
		\label{eq-weight_3_eisenstein}
		q \frac{d}{dq} \log Q = E_{3,m}(\tau) = 1 + \sum_{k=1}^{\infty} k^2 \mathsf{N}_{m,k} \frac{q^k}{1-q^k} \, ,
	\end{equation}
	where $Q$ is the exponentiated local mirror map (\ref{eq-local_mirror_map}) composed with the closed mirror map. It follows that for each $m = 3,4,5,6$, that the integers $\mathsf{N}_{m,k} \in \mathbb{Z}$ are $m$-fold periodic. Similarly, the weight-4 Eisenstein series $E_{4,N}(\tau)$ admit a $q$-Lambert series expansion as
	\begin{equation}
		\label{eq-weight_4_eisenstein}
		q \frac{d}{dq} \log Q =E_{4,N}(\tau) = 1 + \sum_{k=1}^{\infty} k^3 \mathsf{N}_{N,k} \frac{q^k}{1-q^k} \, ,
	\end{equation}
	and that for each $N = 2,3,4,5,6$, that (up to a constant integer multiple $m_N > 0$) integers $m_N \cdot \mathsf{N}_{N,k} \in \mathbb{Z}$ are $N$-fold periodic. We call such numbers \emph{Eisenstein weights}. These integers are listed in Table \ref{table-Eisenstein_weights}. Analogous statements hold for the remaining modular rational elliptic surfaces and pairs $(N,d)$ for Fano threefolds. 
	\renewcommand{\arraystretch}{1.2}
	\begin{table}[!ht]
		\caption{\tiny Periodic Eisenstein weights for the series $E_{3,m}(\tau)$ and $E_{4,N}(\tau)$ from Tables \ref{table-modular_elliptic_data} and \ref{table-K3_modular}. In the latter, weights are listed with the least integer $m_N > 0$.}\label{table-Eisenstein_weights}
		\scalebox{0.8}{
			\begin{tabularx}{.65\textwidth} { 
					>{\hsize=15mm\centering\arraybackslash}X 
					| >{\hsize=75mm\centering\arraybackslash}X  }
				Index & Integral Eisenstein weights   \\
				\hline 
				$m =3$ & $\mathsf{N}_{3,k} =  9,  -9,  0$ \\
				\hline
				$m =4$ &  $\mathsf{N}_{4,k} = 4,  0,  -4, 0$\\
				\hline
				$m =5$& $\mathsf{N}_{5,k} = 2,  1,  -1, -2, 0$ \\
				\hline
				$m =6$ & $\mathsf{N}_{6,k} = 1,  1, 0,  -1,  -1,  0$ \\
				\hline 
				$N =2$ & $1\cdot \mathsf{N}_{2,k} = -80,-40$ \\
				\hline
				$N =3$ & $1\cdot \mathsf{N}_{3,k} =-30,-30,-20$ \\
				\hline 
				$N =4$ & $1 \cdot \mathsf{N}_{4,k} = -16,-16,-16,-12$ \\
				\hline
				$N =5$ & $1 \cdot \mathsf{N}_{5,k} = -10,-10,-10,-10,-8$ \\
				\hline
				$N =6$ & $6 \cdot \mathsf{N}_{6,k}=-42, -39, -44, -39, -42, -34$  \\
				\hline
				$N =\widetilde{6}$ & $6\cdot \mathsf{N}_{\widetilde{6},k} =-48,  -24,  -64,  -24,  -48,  -32$ \\
				\hline
		\end{tabularx}}
	\end{table}
	\par By inverting Equation (\ref{eq-weight_3_eisenstein}) one obtains a functional relationship between the exponentiated mirror maps $Q$ and $q$, determined up to a multiplicative constant $\beta \in \mathbb{Q}^*$. In physics, determining a choice for $\beta$ is called ``fixing the B-field", and this can be quite subtle. In fact, there is a positive integer $\nu \in \mathbb{N}$ such that all quantities can be written in terms of $\mathsf{Q} := Q^\nu$ and $\mathsf{q} := q^\nu$. The functional relationship between $Q$ and $q$ is given by
	\begin{equation}
		\label{eq-functional_rel}
		\beta^{-1}\mathsf{Q} = \mathsf{q} \prod_{k=1}^{\infty} (1- \mathsf{q}^k)^{-k \mathsf{N}_{m,k}} \, , \hspace{10mm} \beta \mathsf{q} = \mathsf{Q} \prod_{k=1}^{\infty} (1- \mathsf{Q}^k)^{-k^2 \widetilde{\mathsf{N}}_{m,k}} \, .
	\end{equation}
	There is an analogous functional relationship obtained by inverting (\ref{eq-weight_4_eisenstein}) and swapping the roles of the numbers in the power of each equation. This leads us to the following. 
	\begin{definition}
		\label{def-weight_4_inverse_expansion}
		For each pair $(N,d)$ in (\ref{eq-degree_level_invariants}) of numerical invariants labeling deformation classes of rank-1 Fano threefolds, define the numbers $\widetilde{\mathsf{N}}_{N,k} \in \mathbb{Q}$ via the Lambert series expansion
		\begin{equation}
			\label{eq-local_fourfold_Q_expansion1}
			\mathsf{Q} \frac{d}{d\mathsf{Q}} \log \mathsf{q} =: 1 + \sum_{k=1}^{\infty} k^2 \widetilde{\mathsf{N}}_{N,k} \frac{\mathsf{Q}^k}{1-\mathsf{Q}^k}
		\end{equation}
		where $\mathsf{Q} = Q^d$ and $\mathsf{q} = q^d$. 
	\end{definition}
	\par Equation (\ref{eq-local_fourfold_Q_expansion1}) is the functional inverse of the weight-4 $\mathbb{Q}$-Eisenstein series $E_{4,N}(\tau)$, pulled back by the morphism $q \mapsto q^d$; in fact, these are reciprocals, as is the case for the analogous weight-3 Eisenstein series. As in (\ref{eq-functional_rel}), the numbers $\widetilde{\mathsf{N}}_{N,k} \in \mathbb{Q}$ are determined up to a choice of the B-field $\beta \in \mathbb{Q}^*$. Moreover, we may normalize (\ref{eq-local_fourfold_Q_expansion1}) by multiplication by some $\mu \in \mathbb{Q}^*$ to obtain an expression like (\ref{eq-local_threefold_Q_expansion1}). 
	\par A priori, the numbers $\widetilde{\mathsf{N}}_{N,k} \in \mathbb{Q}$ are rational, given the rational Eisenstein weights in Table \ref{table-Eisenstein_weights} and the functional relation (\ref{eq-functional_rel}). However, we find for each $(N,d)$ on an experimental basis that up to normalization $\beta,\mu \in \mathbb{Q}^*$, we have $\widetilde{\mathsf{N}}_{N,k} \in \mathbb{Z}$ are actually integral. Thus, we call these numbers \emph{local instanton} numbers, in analogy with the local threefold case. Leading order expansions for the Fano threefolds under consideration here are listed in Table \ref{table-local_instantons}, with $\beta = -1$, $\nu =d= 1$, and $\mu=1$. Analogous computations hold for all other pairs $(N,d)$, whether or not they are naturally related to the twist construction. In fact, the local instanton numbers for $N=2$ appear to agree with known genus zero invariants of local $\mathbb{P}^3$, up to a change in normalization discussed in \S \ref{sss-local_P3}.
	\renewcommand{\arraystretch}{1.2}
	\begin{table}[!ht]
		\caption{\tiny Local instantons for Fano threefolds $(N,1)$ with $N=2,3,4,5,6$. The $N=\widetilde{6}$ is also included. Here $\beta = -1$, $\nu=d=1$ and $\mu =1$. Similar integrality results hold for all $(N,d)$.}\label{table-local_instantons}
		\scalebox{0.8}{
			\begin{tabularx}{.70\textwidth} { 
					>{\hsize=10mm\centering\arraybackslash}X 
					| >{\hsize=95mm\arraybackslash}X  }
				$N$ & \; \; \; \; \; \; \; \; \; Local instantons $\widetilde{\mathsf{N}}_{N,k}$   \\
				\hline 
				$2$ & $-80, \; 3320, \;-272240, \; 29945760, \; -3860155600, \dots$ \\
				\hline
				$3$ & $-30, \; 525, \;-17610, \; 797820, \; -42336600, \dots$ \\
				\hline 
				$4$ & $-16, \; 168, \; -3248, \; 85552, \; -2638480, \dots$ \\
				\hline
				$5$ & $-10, \; 75, \; -980, \; 17640, \; -371690, \dots$ \\
				\hline
				$6$ & $-7,\; 41, \; -400,\; 5425, \; -86150, \; \dots$  \\
				\hline
				$\widetilde{6}$ & $-8, \; 44, \; -448, \; 6288, \; -103192, \dots$ \\
				\hline
		\end{tabularx}}
	\end{table}
	\par Let us briefly discuss how to compute such quantities from our perspective. Consider the canonical local MUM basis of solutions $\omega_0,\omega_1,\omega_2,\omega_3$ at $t=0$ as in (\ref{eq-MUM_basis}) to the quantum differential operator $\mathcal{D}_4 = \theta \cdot \mathcal{L}_3$, obtained from a third order Picard-Fuchs operator $\mathcal{L}_3=\mathcal{L}_{3,N}$. Since $\mathcal{L}_3 = \mathcal{L}_2^{\otimes 2}$ is necessarily a symmetric square, and we have that $\omega_0\omega_2=\omega_1^2$. Writing $\tau = \omega_1/\omega_0$, this quadratic relation is equivalent to $\omega_2=\omega_0 \tau^2$. 
	\par We may always express the new solution $\omega_3$ of $\mathcal{D}_4$ as
	\begin{equation}
		\label{eq-D4_new_sol}
		\omega_3 = \omega_0 \tau^3 + \frac{1}{(2 \pi i)^3} \rho_3 \, ,
	\end{equation}
	where $\rho_3(t)$ is a power series defined by 
	\begin{equation}
		\label{eq-rho3_def}
		\rho_3(t) = h_3(t) - \frac{h_1(t)^3}{h_0(t)^2} \, .
	\end{equation}
	Thus $\rho_3$ measures the failure of $\mathcal{D}_4$ to be the symmetric cube of a second order operator. Then the rank-3 Eisenstein potential (\ref{eq-eisenstein_potential}) is $\mathscr{E}_{3,N} = \frac{\omega_3}{\omega_0}$, and the associated weight-4 Eisenstein series is given by 
	\begin{equation*}
		E_{4,N}(\tau) = (2\pi i)^3 \frac{d^3}{d\tau^3}\mathscr{E}_{3,N}(\tau) \, .
	\end{equation*}
	Thus one only needs the additional piece of data $h_3(t)$. Such functions have arisen e.g. in Golyshev and Zagier's proof of the Gamma conjecture for rank-1 Fano threefolds \cite[\S 2.1]{MR3462676}. For convenience, we list these in Table \ref{table-frobenius_solutions}. An analogous story holds in the local threefold case utilizing Theorem \ref{thm-eisenstein}. 
	\renewcommand{\arraystretch}{1.4}
	\begin{table}[!ht]
		\caption{\tiny Leading order expansion of canonical frobenius solution $h_3(t)$ from the fourth order operator $\mathcal{D}_{4,N}$, $N=2,3,4,5,6,\widetilde{6}$.}\label{table-frobenius_solutions}
		\scalebox{.8}{
			\begin{tabularx}{.70\textwidth} { 
					>{\hsize=10mm\centering\arraybackslash}X 
					| >{\hsize=95mm\arraybackslash}X  }
				$N$ & \; \; \; \; \; \; \; \; \; Frobenius solution $h_3(t)$   \\
				\hline 
				$2$ & $ -480 t-61740 t^2-\frac{77788384 t^3}{9}-\frac{8337006209 t^4}{6}- \cdots$ \\
				\hline
				$3$ & $-180 t-\frac{19845 t^2}{2}-\frac{1775576 t^3}{3}-\frac{644909189 t^4}{16}- \cdots$ \\
				\hline 
				$4$ & $-96 t-3180 t^2-\frac{1019552 t^3}{9}-\frac{82535939 t^4}{18}- \cdots$ \\
				\hline
				$5$ & $-60 t-\frac{2775 t^2}{2}-\frac{310781 t^3}{9}-\frac{15547863 t^4}{16}- \cdots$ \\
				\hline
				$6$ & $-42 t-\frac{3033 t^2}{4}-\frac{522389 t^3}{36}-\frac{90124865 t^4}{288}- \cdots$  \\
				\hline
				$\widetilde{6}$ & $-48 t-990 t^2-\frac{179044 t^3}{9}-\frac{5395165 t^4}{12}- \cdots$ \\
				\hline
		\end{tabularx}}
	\end{table}
	\subsection{Computing local invariants}
	\label{ss-local_invariants}
	\par As discussed in \S \ref{ss-Fano_MS_&_QDE}, the mirrors of anticanonical elliptic curves in $\mathbb{P}^2$ and anticanonical quartic K3 surfaces in $\mathbb{P}^3$ can be connected via the twist construction with generalized functional invariant $(i,j,\alpha)=(3,1,1)$. 
	\subsubsection{Local \texorpdfstring{$\mathbb{P}^2$}{P^2}}
	\label{sss-local_P2}
	The inversion scheme of (\ref{eq-functional_rel}) applied to the modular form $E_{3,3}(\tau)$ with $\beta = -1$, $\mu = \frac{1}{3}$, and $\nu = 3$ yields the instanton numbers (here $\widetilde{\mathsf{N}}_k \equiv \frac{1}{3} \widetilde{\mathsf{N}}_{3,k}$)
	\begin{equation}
		\label{eq-KP2_GV0_invariants}
		\scalemath{0.7}{
			\begin{array}{lcl}
				\widetilde{\mathsf{N}}_1=3, & \hspace{15mm} & \widetilde{\mathsf{N}}_9= 27748899, \\
				\widetilde{\mathsf{N}}_2=-6, & \hspace{15mm} & \widetilde{\mathsf{N}}_{10}=-360012150, \\
				\widetilde{\mathsf{N}}_3=27, & \hspace{15mm} & \widetilde{\mathsf{N}}_{11}=4827935937, \\
				\widetilde{\mathsf{N}}_4=-192, & \hspace{15mm} & \widetilde{\mathsf{N}}_{12}=-66537713520, \\
				\widetilde{\mathsf{N}}_5= 1695, & \hspace{15mm} & \widetilde{\mathsf{N}}_{13}=938273463465, \\
				\widetilde{\mathsf{N}}_6=-17064, & \hspace{15mm} & \widetilde{\mathsf{N}}_{14}=-13491638200194, \\
				\widetilde{\mathsf{N}}_7= 188454, & \hspace{15mm} & \widetilde{\mathsf{N}}_{15}=197287568723655, \\
				\widetilde{\mathsf{N}}_8=-2228160, & \hspace{15mm} & \;\;\;\;\;\;\;\vdots
		\end{array}}
	\end{equation}
	in perfect agreement with the literature, such as Coates and Iritani \cite[Table 2]{MR4301054}, for the genus zero GV invariants of $K_{\mathbb{P}^2}$. That this identification persists for all $k \in \mathbb{N}$ is a consequence of Bloch, Kerr, and Vanhove \cite{MR3780269}. Moreover, with alternative rational normalizations such as  $\mu = 3$, one computes the instanton expansion for the local Calabi-Yau threefold $K_{B_6}$ \cite[Eqn. (4.10)]{MR1468700}, where it is referred to as the so-called $E_6$ noncritical string. Similarly, local instantons for $K_{B_5}$ are computed from the modular form $E_{3,4}(\tau)$ with $\beta=-1$, $\mu =4$, and $\nu = 2$. The so-called conifold \cite[Eqn. (4.10)]{MR1468700} is obtained from the same datum but with $\beta=1$. It is expected, as with Stienstra \cite{MR2282958} and the results in loc. cit., that such an identification persists for all $k \in \mathbb{N}$; note that the underlying del Pezzo surface in these examples is in general not toric. 
	\subsubsection{Local \texorpdfstring{$\mathbb{P}^3$}{P^3}}
	\label{sss-local_P3}
	The case of the rank-1 Fano threefold $\mathbb{P}^3$, whose deformation space is obviously trivial, corresponds to the pair $(N, d)=(2,4)$. The relevant Picard-Fuchs operator is the pullback of $\mathcal{L}_{3,2}$ in Table \ref{table-modular_pf} under the map $t \mapsto t^4$ which induces a branched covering of the coarse moduli space $X_0(2)^+$ of $\mathrm{M}_2$-polarized K3s; see \cite[\S 5.8]{MR2306141}. From Definition (\ref{def-weight_4_inverse_expansion}), we take $\nu = 4$, and after computing the relevant expansions from Equation (\ref{eq-local_fourfold_Q_expansion1}), one determines that by choosing $\beta =1$ and $\mu =-\frac{1}{4}$, the following normalized local instanton numbers (here $\widetilde{\mathsf{N}}_k \equiv -\frac{1}{4}\widetilde{\mathsf{N}}_{4,k}$) are obtained:
	\begin{equation}
		\label{eq-KP3_GW0_invariants}
		\scalemath{0.7}{
			\begin{array}{lcl}
				\widetilde{\mathsf{N}}_1=-20, & \hspace{15mm} & \widetilde{\mathsf{N}}_7=-21025364147340,  \\
				\widetilde{\mathsf{N}}_2=-820, & \hspace{15mm} & \widetilde{\mathsf{N}}_8=-3381701440136400, \\
				\widetilde{\mathsf{N}}_3=-68060, & \hspace{15mm} & \widetilde{\mathsf{N}}_9=-565563880222390140, \\
				\widetilde{\mathsf{N}}_4=-7486440, & \hspace{15mm} & \widetilde{\mathsf{N}}_{10}=-97547208266548098900, \\
				\widetilde{\mathsf{N}}_5=-965038900, & \hspace{15mm} & 	\widetilde{\mathsf{N}}_{11}=-17249904137787210605980,\\
				\widetilde{\mathsf{N}}_6=-137569841980, & \hspace{15mm} & \widetilde{\mathsf{N}}_{12}=-3113965536138337597215480,  \\
				& \hspace{15mm} & \;\;\;\;\;\;\; \vdots \\
		\end{array}}
	\end{equation}
	Then for each $k=1,\dots,12$, the integers $\widetilde{\mathsf{N}}_k$ agree precisely with the genus zero invariants for the local Calabi-Yau fourfold $K_{\mathbb{P}^3}$, the total space of the canonical bundle of $\mathbb{P}^3$. These were computed for example by Klemm and Pandharipande \cite[Tbl. 1]{MR2415462} using the Aspinwall-Morrison formula \cite{MR1204770}, and the bound of $k=12$ comes from loc. cit., the most extensive list we could find in the literature relevant to this case. We expect that this agreement holds for $k \in \mathbb{N}$.
	\subsection{A conjecture on local invariants}
	\label{ss-conj}
	\par Based on the experimental results in \S \ref{sss-local_P3}, we make the following conjecture.
	\vspace{3mm}
	\begin{flushleft}
		\textbf{Conjecture.} (a) For each deformation class of rank-1 Fano threefold $F$ of type $(N, d)$ in Equation (\ref{eq-degree_level_invariants}), there are rational numbers $\beta,\mu \in \mathbb{Q}^*$ such that the normalized local instanton numbers $\mu\widetilde{\mathsf{N}}_{N,k}$ from Equation (\ref{eq-local_fourfold_Q_expansion1}) expressed in $\mathsf{Q} = Q^d$ with choice of B-field determined by $\beta$ are genus zero invariants of the local Calabi-Yau fourfold $K_{F}$. Moreover, (b) the $\mathsf{Q}$-expansions in Equation (\ref{eq-local_fourfold_Q_expansion1}) are expressible in terms of quasi-modular forms for some discrete subgroup of $\mathrm{SL}(2,\mathbb{R})$ commensurable with $\mathrm{SL}(2,\mathbb{Z})$.
	\end{flushleft}
	\begin{remark}[Quasi-modularity, The local-relative correspondence]
		\label{rem-local_relative}
		Quasi-modularity of instanton expansions in string theory, especially within the framework of local mirror symmetry, has been observed for local Calabi-Yau threefolds by Alim, Scheidegger, Yau, and Zhou \cite{MR3273318}, and also for Calabi-Yau fourfolds in \cite{MR3733768}. The local-relative correspondence of van Garrel, Graber, and Ruddat \cite{MR3948687}, together with the work on relative invariants by Fan, Tseng, and You \cite{MR3997137}, and You ~\cite{MR4716739} would then imply that our computations of appropriately normalized local instantons also compute relative invariants, up to a sign.
	\end{remark}
	It is expected that such an analogous observation for the case $n=2$ is consistent with a result of Doran and Kerr \cite[Cor. 5.3]{MR2851155} on the asymptotic growth rate of local GW invariants.

	\newpage
	
	\bibliographystyle{abbrv}
	\bibliography{MalmendierSchultz_bib}
	
	\appendix

	\section{Modular differential equations}
	\label{s-modular_DEs}
	\par Let $\mathbb{H} \subset \mathbb{C}$ be the upper half plane and $\mathbb{H}^* =  \mathbb{H} \cup \mathbb{Q}  \cup \{i\infty\}$. We use standard notations for the congruence subgroups $\Gamma_i(N) \subset \mathrm{SL}(2,\mathbb{Z})$ and the associated modular curves $X_j(N) = \Gamma_j(N) \backslash \mathbb{H} \cup \{\mathrm{cusps}\}$, $j=0,1$. Similarly $M_k(\Gamma)$ will represent the space of weight-$k$ modular forms associated with the  subgroup $\Gamma$. Let $w_N \in \mathrm{SL}(2,\mathbb{R})$ be the Fricke involution mapping $\tau \in \mathbb{H}$ to $-1/N\tau \in ~\mathbb{H}$. We obtain the augmented group $\Gamma_0(N)^+ \subset \mathrm{SL}(2,\mathbb{R})$ generated by $\Gamma_0(N) \cup \{w_N\}$, and the corresponding modular curve $X_0(N)^+ = \Gamma_0(N)^+ \backslash \mathbb{H} \cup \{\mathrm{cusps}\}$. Then $\Gamma_0(N)^+\subset \mathrm{SL}(2,\mathbb{R})$ is a discrete subgroup commensurable with the modular group $\mathrm{SL}(2,\mathbb{Z})$. The modular curves $X_0(N)$ and $X_0(N)^+$ respectively parameterize pairs of ordered and unordered pairs of elliptic curves related by cyclic degree-$N$ isogeny. The modular curve $X_1(N)$ parameterizes elliptic curves with a rational $N$-torsion point. There only finitely many values of $N$ for which the modular curves $X_0(N),\, X_0(N)^+,$ and $X_1(N)$ are genus zero; see, for example, Maier \cite{MR2514149}.
	\par Suppose that $X = \Gamma \backslash \mathbb{H} \cup \{\mathrm{cusps}\} \cong \mathbb{P}^1$ is such a genus-zero modular curve, where $\Gamma \subseteq \mathrm{SL}(2,\mathbb{R})$ is a discrete subgroup that is commensurable with $\mathrm{SL}(2,\mathbb{Z})$. Being rational, the function field $\mathbb{C}(X) \cong \mathbb{C}(t)$ is isomorphic to the field of rational functions in a local single valued algebraic coordinate $t \in X$; such a local coordinate is called a \emph{Hauptmodul} of the modular curve $X$. Since $X$ is a modular curve, $t$ is itself a modular function, i.e., a weight-0 modular form. Let $\pi\colon \mathbb{H} \to X$ be the canonical projection. A multivalued inverse map $\tau \colon X \to \mathbb{H}$ is called a \emph{uniformization} of $X$. 
	\par The multivalued function $\tau(t)$ is studied by the classical Gauss-Schwarz theory. In \cite{MR0429918}, Shioda shows -- following results of Kodaira \cite{MR0165541} -- that the modular curve $X$ naturally serves as the base of a pencil $\pi : \mathcal{E} \to X$ of elliptic curves, where the fibre $\pi^{-1}(t) = E_t$ is an elliptic curve over $\mathbb{C}$ is constructed from the relevant modular subgroup $\Gamma$; specific details on elliptic fibrations are provided in the Appendix \ref{s-twist}. This can be seen as follows. Let $\mathcal{J} = \mathrm{SL}(2,\mathbb{Z}) \backslash \mathbb{H} \cup \{\mathrm{cusps}\} \cong \mathbb{P}^1$ be the $j$-line, whose Hauptmodul is Klein's modular $j$-function
	\begin{equation}
		\label{eq-j_invariant}
		j(\tau) = 1728\frac{E_4(\tau)^3}{E_4(\tau)^3-E_6(\tau)^2} \, ,
	\end{equation}
	where $E_4,E_6$ are the Eisenstein series
	\begin{equation}
		\label{eq-eisenstein_series_Ek}
		E_k(\tau) = 1 - \frac{2k}{B_k} \sum_{n=1}^{\infty} \frac{n^{k-1} q^n}{1-q^n} \; ,
	\end{equation}
	for $k=4,6$ and $q = \exp(2\pi i \tau)$. Here, $B_k$ is the $k^{th}$ Bernoulli number. The series are holomorphic modular forms of weight-4 and weight-six, respectively, for the modular group $\mathrm{SL}(2,\mathbb{Z})$, and the quantity $E_2(\tau)$ defined by the same formula for $k=2$ defines a quasi modular form. We will also use the usual Dedekind $\eta$-function
	\begin{equation}
		\label{eq-dedekind_eta}
		\eta(\tau) = q^{\frac{1}{24}} \prod_{n=1}^{\infty}(1-q^n) \; .
	\end{equation}
	\par There is a canonical map $J \colon X \to \mathcal{J}$ coming from the inclusion $\Gamma \hookrightarrow \mathrm{SL}(2,\mathbb{Z})$ and quotient spaces, resulting from the natural action on the upper half plane $\mathbb{H}$. Let $U = X - \{ \mathrm{cusps}\}$ and fix $t_0 \in U$, letting $G = \pi_1(U,t_0)$. Let $\pi^\prime \colon U^\prime \to U$ be the universal cover. Then there is a holomorphic map $\widetilde{\omega} \colon U^\prime \to \mathbb{H}$ making the following diagram commutative,
	\begin{equation}
		\label{eq-functional_invariant_diagram}
		\begin{tikzcd}
			U^\prime \arrow{r}{\widetilde{\omega}} \arrow{d}[left]{\pi^\prime} & \mathbb{H} \arrow{d}{\pi} \arrow{r}{j} & \mathcal{J} \arrow[equal]{d} \\
			U \arrow[hookrightarrow]{r}{\iota} & X \arrow{r}{J} &  \mathcal{J}
		\end{tikzcd}
	\end{equation}
	where $\iota \colon U \hookrightarrow X$ is the inclusion map. Moreover, there is a unique representation 
	\begin{equation*}
		\varphi \colon \quad G \to \Gamma \subseteq \mathrm{SL}(2,\mathbb{Z})
	\end{equation*}
	lifting the action of the fundamental group $G$ as deck transformations on the universal covering $\pi^\prime \colon U^\prime \to U$ to the natural action of $\Gamma$ on $\mathbb{H}$, i.e., $\omega(g \cdot u^\prime) = \varphi(g) \cdot \omega(u^\prime)$ for all $g \in G$ and $u^\prime \in U^\prime$. Hence, the representation $\varphi$ furnishes $X$ with a sheaf $\mathcal{G} \to X$ that is locally constant over $U$, with generic stalk given by $\mathbb{Z} \oplus \mathbb{Z}$ on which $\Gamma$ acts by monodromy representations as the transition functions. 
	\par Results by Kodaira \cite{MR0165541} guarantee that there is an elliptic pencil $\mathcal{E} \to X$ realizing the data above in the following way: the pencil is a locally topologically trivial fibration whose generic fibre $E_t$ is a smooth elliptic curve, such that the $j$-invariant of $E_t$ coincides with $J(t) \in \mathbb{P}^1$, and $H_1(E_t,\mathbb{Z}) \cong \mathcal{G}_t = \mathbb{Z} \oplus \mathbb{Z}$ generates the sheaf $\mathcal{G} \to X$. The elliptic pencil can be represented concretely as a \emph{Weierstrass model}, where the fibre $E_t$ is represented by the affine cubic
	\begin{equation}
		\label{eq-weierstrass}
		E_t \colon \quad y^2 = 4x^3 -g_2(t)x - g_3(t) \, ,
	\end{equation}
	where $g_2(t), \, g_3(t)$ are sections of a certain line bundle described in more detail in Appendix \ref{ss-wm}. The discriminant of the right hand side is $\Delta(t) = g_2(t)^3 - 27g_3(t)^2$, determining where $E_t$ becomes singular. This discriminant locus contains the cusps of the base curve. Hence, the modular subgroup $\Gamma$ acts on $H_1(E_t,\mathbb{Z})$ by monodromy transformations as one moves around the discriminant locus. Thus Kodaira's result says that the sheaf $\mathcal{G}$ and holomorphic map $J \colon X \to \mathcal{J} \cong \mathbb{P}^1$ are of geometric origin, they are respectively the \emph{homological invariant} and \emph{functional invariant} of a family of elliptic curves parameterized by $X$. 
	\par We may assume without loss of generality that $t=0$ is a singular point. For $t \in \mathbb{P}^1$ sufficiently close to zero, we can choose a locally constant basis $\Sigma_0,\Sigma_1 \in H_1(E_t,\mathbb{Z})$ of integral 1-cycles for the homological invariant $\mathcal{G}$ on which a fiberwise non-vanishing holomorphic 1-form $\Omega_t \in H^{1,0}(E_t)$ can be integrated, yielding the period integrals 
	\begin{equation*}
		\omega_i(t) = \int_{\Sigma_i} \Omega_t
	\end{equation*}
	for $i = 0, 1$, as in Equation (\ref{eq-periods}). Such a choice of $\Omega_t \in H^{1,0}(E_t)$ can be made to vary holomorphically in $t$, and is called an \emph{analytic marking} of the pencil $\mathcal{E} \to X$. Then the period integrals become generically multivalued holomorphic functions of $t$ and satisfy a second order Picard-Fuchs equation
	\begin{equation}
		\label{eq-second_order_PF}
		\mathcal{L}_2(\omega)= a_2(t)\frac{d^2 \omega}{dt^2} + a_1(t)\frac{d \omega}{dt} + a_0(t)\omega =0 \; .
	\end{equation}
	The Picard-Fuchs operator $\mathcal{L}_2$ can be computed from the geometric data in the Weierstrass model in Equation \ref{eq-weierstrass}. It follows from well known results (see, for example, Stiller \cite[\S 3]{MR0607940}) that the periods $\omega_0(t), \, \omega_1(t)$ satisfy the differential relations
	\begin{equation}
		\label{eq-elliptic_pf}
		\renewcommand{\arraystretch}{1.3}
		\frac{d}{d t}\left(\begin{array}{c}
			\omega_0 \\
			\omega_1
		\end{array}\right)=
		\left(\begin{array}{cc}
			-\frac{1}{12} \frac{\Delta^{\prime}}{\Delta} & \frac{3}{2}\frac{ \delta}{ \Delta} \\
			-\frac{1}{8}\frac{g_2 \delta}{ \Delta} & \frac{1}{12} \frac{\Delta^{\prime}}{\Delta}
		\end{array}\right) \cdot\left(\begin{array}{c}
			\omega_0 \\
			\omega_1
		\end{array}\right) \, ,
	\end{equation}
	where $\Delta(t)$ is the discriminant and $\delta(t) = 3g_3(t)g_2(t)^{\prime} - 2g_2(t)g_3(t)^{\prime}$, which is equivalent to a second order equation, i.e.,  (\ref{eq-second_order_PF}). Equation (\ref{eq-elliptic_pf}) is the Gauss-Manin connection of the elliptic pencil. Moreover, the monodromy group of this equation is precisely $\Gamma$, see, for example, Yang and Zudilin \cite[\S 1]{MR2731075} for a discussion relevant to the contents of this article. We can always choose an ordering of the periods $\omega_0,\omega_1$ (after relabeling if necessary) such that the image 
	\begin{equation*}
		\omega(U) \subseteq \mathbb{H}^* \subset \mathbb{P}^1 \; ,
	\end{equation*} 
	\begin{equation}
		\label{eq-projective_period_map}
		\omega(t) = \left[\omega_0(t) \, : \, \omega_1(t)\right] = \left[1 \, : \frac{\omega_1(t)}{\omega_0(t)}\right] \equiv \left[1 \, : \, \tau(t)\right] \; ,
	\end{equation}
	such that the image of $t=0 \in X$ is mapped to $i\infty \in \mathbb{H}^* \subset \mathbb{P}^1$. The resulting holomorphic map $\tau \colon U \to \mathbb{H}$ in Equation (\ref{eq-projective_period_map}) is called the \emph{projective} period map of the modular elliptic pencil $\mathcal{E} \to X$, and gives a multivalued inverse to the covering map $\pi^{\prime} \colon U^{\prime} \to U$ in Equation (\ref{eq-functional_invariant_diagram}), in particular $\tau(U) = \widetilde{\omega}(U^{\prime}) \subseteq \mathbb{H}$. This implies that $J(t) = j(\tau(t))$ for all $t \in U$, in congruence with the discussion above.
	\par Assuming that $\Gamma$ contains the $\mathbb{H}$-automorphism $\tau \mapsto \tau + 1$, we can express the Hauptmodul $t$ in a $q$-series about $i \infty \in \mathbb{H}^*$. One obtains
	\begin{equation}
		\label{eq-Hauptmodul_q_expansion}
		t=t(\tau) = \sum_{k=1}^\infty c_k q^k
	\end{equation}
	by inverting the map $\tau(t) = \omega_1(t)/\omega_0(t)$ term by term; this is precisely the mirror map of the elliptic pencil $\mathcal{E} \to X$ from Equation (\ref{eq-mirror_map}). Under favorable circumstances, in a neighborhood of $t=0$, the period $\omega_0(t)$ is a single-valued holomorphic function, and expressing $\omega_0=\omega_0(\tau)$ via the $q$-series in Equation (\ref{eq-Hauptmodul_q_expansion}) yields a $\Gamma$-modular form of weight one. We then say that the Picard-Fuchs equation in Equation (\ref{eq-second_order_PF}) admits a \emph{modular} parameterization. Let $F(\tau) = \omega_0(t(\tau)) \in M_1(\Gamma)$ be the modular period. Since 
	\begin{equation*}
		\omega_1 = \tau \omega_0 = \tau F(\tau) 
	\end{equation*}
	is also a solution of the Picard-Fuchs equation that is linearly independent from $F(\tau)$, it follows that under the change of variables in Equation (\ref{eq-Hauptmodul_q_expansion}), $\mathcal{L}_2$ transforms into the particularly simple second order operator
	\begin{equation}
		\label{eq-modular_PF}
		\mathcal{L}_2 = \frac{1}{E_3(\tau)}\left(\frac{1}{2\pi i} \, \frac{d}{d\tau}\right)^2 \frac{1}{F(\tau)} \, .
	\end{equation}
	Here, a straightforward computation shows that $E_3(\tau)$ is the product
	\begin{equation}
		\label{eq-yukawa_modular_form}
		E_3(\tau) = F(\tau) \frac{t^{\prime}(\tau)}{t(\tau)}  \, ,
	\end{equation}
	where $t^{\prime} = \frac{1}{2\pi i}\frac{dt}{d\tau}$. Importantly for this article, $E_3(\tau) \in M_3(\Gamma)$ transforms as a  $\Gamma$-modular form of weight three, in particular, it is precisely the Eisenstein series from (\ref{eq-eisenstein_factorization}).
	\par More generally, we have the following based on a folklore result going back to Gauss, Fricke, Klein, etc. For a related modern discussion, see Zagier \cite[\S 5.4]{MR2409678}, and also Yang and Zudilin \cite[Theorem B]{MR2731075}.
	\begin{theorem}
		\label{thm-modular_ODEs}
		Let $X = \Gamma \backslash \mathbb{H} \cup \{\mathrm{cusps}\}$ be a genus-zero modular curve for a discrete subgroup $\Gamma \subseteq \mathrm{SL}(2,\mathbb{R})$ commensurable with $\mathrm{SL}(2,\mathbb{Z})$, and let $t \in \mathbb{C}(X)$ be a Hauptmodul. Suppose that 
		\begin{equation}
			\label{eq-fuchsian_ODE}
			\mathcal{L}_n = a_n(t)\frac{d^n}{dt^n}+a_{n-1}(t)\,\frac{d^{n-1}}{dt^{n-1}} + \cdots + a_1(t)\,\frac{d}{dt}+a_0(t) 
		\end{equation}
		is a Fuchsian operator of order $n=2$ or $n=3$ defined on $X - \{\mathrm{cusps}\}$, and $t=0$ is a singular point of maximally unipotent monodromy, with monodromy group determined by the natural 2-dimensional representation of $\Gamma$ in $\mathrm{SL}(2,\mathbb{R}) \cong \mathrm{Sp}(2,\mathbb{R})$ if $n=2$, or as a symmetric square representation in $\mathrm{SO}(2,1)$, if $n=3$. Let $\omega_0(t),\omega_1(t)$ be the holomorphic and logarithmic solutions at $t=0$ as in Equation (\ref{eq-MUM_basis}), and $\tau(t)=\omega_1(t)/\omega_0(t)$. If the $q$-expansion of the holomorphic solution $\omega_0(\tau)=F(\tau) \in M_{n-1}(\Gamma)$ at $t=0$ is a $\Gamma$-modular form of weight $(n-1)$, then under change of variables $t=t(\tau)$ as in Equation (\ref{eq-Hauptmodul_q_expansion}), $\mathcal{L}_n$ is projectively equivalent to the operator
		\begin{equation}
			\label{eq-modular_CY_op}
			\mathcal{L}_n = \frac{1}{E_{n+1}(\tau)}\left(\frac{1}{2\pi i} \, \frac{d}{d\tau}\right)^n \frac{1}{F(\tau)} \, ,
		\end{equation}
		where $E_{n+1}(\tau) = F(\tau) t^{\prime}/t \in M_{n+1}(\Gamma)$ is a weight-$(n+1)$-modular form and $t^{\prime} = \frac{1}{2\pi i}\frac{dt}{d\tau}$. This modular form is precisely the Eisenstein series appearing in Theorem \ref{thm-eisenstein}. A fundamental basis of solutions to Equation (\ref{eq-modular_CY_op}) is
		\begin{equation}
			F(\tau), \, \tau F(\tau), \, \dots, \, \tau^{n-1}F(\tau) .
		\end{equation}
	\end{theorem}
	In fact, the converse of this theorem holds, but we do not need it in the context of this article. 
	\section{Weierstrass models and the twist construction}
	\label{s-twist}
	The twist construction of Doran and Malmendier  \cite{MR4069107} provides a geometric mechanism for constructing CY varieties iteratively that admit fibrations useful in string dualities, while also realizing monodromy tuples determined by Bogner and Reiter \cite{MR2980467} that are associated to certain subfamilies of known CY operators  \cite{almkvist_tables_2010}. More specifically, they show how Doran's generalized functional invariant \cite{MR1779161} (defined in Appendix \ref{ss-twisted_wm}) can be used in a ``bottom-up" fashion to iteratively construct families of $m$-dimensional CY varieties $V^{(m)}$ that are simultaneously elliptically fibered over a rational base and fibered by the family of $(m-1)$-dimensional CY varieties $V^{(m-1)}$ over a punctured $\mathbb{P}^1$ such that $t=0$ is a MUM point for each family. Additionally, their construction naturally shows that the Picard-Fuchs operator $\mathcal{L}_n$ for $V^{(m)}$ is given by the Hadamard product of the Picard-Fuchs operator $\mathcal{L}_{n-1}$ for $V^{(m-1)}$ and a certain hypergeometric operator determined by the generalized functional invariant.
	\subsection{Weierstrass models}
	\label{ss-wm}
	\par Central to the twist construction are elliptic fibrations and associated Weierstrass models. If $V$, $S$ are normal complex algebraic varieties of respective dimensions $m$ and $m-1$, an elliptic fibration of $V$ over $S$ is a proper surjective morphism $\pi \colon V \to S$ with connected, irreducible fibers of arithmetic genus one, i.e., each fibre is either a smooth elliptic curve, a rational curve with a node, or a rational curve with a cusp. We assume that $\pi$ is smooth over a Zariski open subset $S_0 \subset S$ such that the complement $D_{\pi} = S - S_0$ is a divisor with at worst normal crossings, and that the nodal or cuspidal fibers occur over $D_{\pi}$. Hence, the generic fibre of $\pi$ is a smooth elliptic curve. Moreover, we assume the fibration is a \emph{Jacobian} elliptic fibration, that is, we assume existence of a section $\sigma \colon S \to V$ marking the neutral point of the elliptic curve group law on each smooth fibre. The fibration is assumed to be minimal in the following way in that there are no $(-1)$-curves in the fibers of $\pi$. A more general notion holds for higher dimensional elliptic fibrations that we will not need here.
	\par Many geometric aspects of elliptic fibrations are determined by the fundamental line bundle $\mathsf{L} := (R^1\pi_*\mathcal{O}_V)^{-1} \in \mathrm{Pic}(S)$, which is related to the situation above to the canonical bundles $K_V,K_S$ of $V$ and $S$ by the adjunction formula as 
	\begin{equation}
		\label{eq-canonical_bundle}
		K_V \cong \pi^*(K_S \otimes \mathsf{L}) \, .
	\end{equation} 
	Thus, the canonical bundle $K_V \cong \mathcal{O}_V$ is trivial when $\mathsf{L} \cong K_S^{-1}$ is the anticanonical bundle of $S$, giving us a way to determine necessary conditions for when the total space $V$ of the fibration is a CY variety. Additionally, the fundamental line bundle determines a Weierstrass model $\hat{\pi} \colon W \to S$ for the fibration as follows. 
	\par Given an elliptic fibration $\pi \colon V \to S$ with fundamental line bundle $\mathsf{L}$, let $g_2, \, g_3$ be sections of $\mathsf{L}^{\otimes 4}, \, \mathsf{L}^{\otimes 6}$ such that the discriminant section $\Delta = g_2^3-27g_3^2$ of $\mathsf{L}^{\otimes 12}$ does not vanish identically. The Weierstrass model $W$ for the fibration is defined as the subvariety 
	\begin{equation}
		\label{eq-Weierstrass_model}
		y^2z=4x^3 - g_2 x z^2-g_3z^3
	\end{equation}
	of the $\mathbb{P}^2$-bundle $\mathbf{P} = \mathbb{P}(\mathcal{O}_S \oplus \mathsf{L}^{\otimes 2} \oplus \mathsf{L}^{\otimes 3})$ over $S$ in the following way. Let $p \colon \mathbf{P} \to S$ be the natural projection, and $\mathcal{O}_{\mathbf{P}}(1)$ be the tautological line bundle. The variables $x, y, z$ in Equation (\ref{eq-Weierstrass_model}) are respectively sections of $\mathcal{O}_{\mathbf{P}}(1) \otimes \mathsf{L}^{\otimes 2}, \, \mathcal{O}_{\mathbf{P}}(1) \otimes \mathsf{L}^{\otimes 3}$, and $\mathcal{O}_{\mathbf{P}}(1)$ that correspond to the natural injections of $\mathsf{L}^{\otimes 2}, \, \mathsf{L}^{\otimes 3}$, and $\mathcal{O}_S$ into $p_*\mathcal{O}_{\mathbf{P}}(1) = \mathcal{O}_S \oplus \mathsf{L}^{\otimes 2} \oplus \mathsf{L}^{\otimes 3}$. Then $\hat{\pi} \colon W \to S$ is naturally an elliptic fibration over $S$ in the sense outlined above, where $\hat{\pi}$ is the natural projection, and the discriminant divisor $D_{\hat{\pi}}$ is given by $\mathcal{O}_S(\Delta)$. The Weierstrass model naturally is a Jacobian elliptic fibration, with canonical section $\hat{\sigma} \colon S \to W$ given by the point $[x \, : \, y \, : \, z] = [0 \, : \, 1 \, : \, 0]$ in each fiber. Moreover, the image $\Sigma := \hat{\sigma}(S) \subset W \subset \mathbf{P}$ is a Cartier divisor, whose normal bundle is isomorphic to the fundamental line bundle as $\hat{\pi}_*\mathcal{O}_{\mathbf{P}}(-\Sigma) \cong \mathsf{L}$.
	\par However, the Weierstrass model $W$ may be singular. Miranda \cite{MR0690264} studied the resolution of the singularities of $W$ in the case that $V$ is a threefold, and shows there is a natural birational map $V \dashrightarrow W$ compatible with the sections $\sigma, \hat{\sigma}$ of the respective fibrations by blowing down all components of the fibers that do not meet $\sigma(S)$. If $\pi \colon V \to S$ is minimal in the sense described above, then the inverse map $W \dashrightarrow V$ is a resolution of the singularities of $W$. In this way, all relevant geometry of the Jacobian elliptic fibration is captured by its associated Weierstrass model.
	\par The twist construction of Doran and Malmendier  is concerned with 1-parameter \emph{families} of Weierstrass models $\mathcal{W} \to \mathbb{P}^1$, where each $t \in \mathbb{P}^1$ corresponds to a Weierstrass model $\hat{\pi} \colon W_t \to S$ for a 1-parameter family $\mathcal{V} \to \mathbb{P}^1$ of $m$-dimensional elliptic CY varieties $\pi \colon V_t \to S$. Then taking local coordinates $u = (u_1,\dots, u_{m-1}) \in S$ and $t \in \mathbb{P}^1$, the Weierstrass model in Equation (\ref{eq-Weierstrass_model}) for $W_t$ becomes
	\begin{equation}
		\label{eq-Weierstrass_model_family}
		y^2=4x^3-g_2(u, t)x-g_3(u, t) \, ,
	\end{equation}
	where we have taken the affine chart $z=1$ in $W_t$. In addition, a holomorphic $m$-form $\Omega_t \in H^{m,0}(V_t)$ can be represented on the Weierstrass model by the closed differential form 
	\begin{equation}
		\label{eq-holomorphic_form_WM}
		\Omega_t = du_1 \wedge \cdots \wedge du_{m-1} \wedge \frac{dx}{y} \, ,
	\end{equation}
	where we suppress the pullback by the birational map $V \dashrightarrow W$. As described in \S \ref{sss-VHS}, this form can then be integrated over a suitable basis of $m$-cycles to obtain the period integrals $\omega_i(t)$ that satisfy the Picard-Fuchs equation as a consequence of Griffiths transversality. See \cite[\S 5.3]{MR4069107} and \cite[\S 2.1]{MR4494119} for more details. As a result of this discussion, we will blur the lines between $V_t$ and the Weierstrass model $W_t$ and hence the families $\mathcal{V} \to \mathbb{P}^1$ and $\mathcal{W} \to \mathbb{P}^1$ through the duration of this article.
	\subsection{Twisted families of Weierstrass models}
	\label{ss-twisted_wm}
	\par We first recall Doran's \emph{generalized functional invariant} \cite{MR1877754}. A generalized functional invariant is a triple $(i, j, \alpha)$ with $i, j \in \mathbb{N}$ and $\alpha \in \left \{\frac{1}{2},1 \right \}$ such that $1 \leq i, j \leq 6$. The generalized functional invariant encodes a 1-parameter family of degree $i + j$ covering maps $\mathbb{P}^1 \to \mathbb{P}^1$,  which is totally ramified over $0$, ramified to degrees $i$ and $j$ over $\infty$, and simply ramified over another point $\mathsf{t}$.  For homogeneous coordinates $[v_0 : v_1] \in \mathbb{P}^1$, the relevant family of maps parameterized by $\mathsf{t} \in \mathbb{P}^1 - \{0,1,\infty\}$ is given by 
	\begin{equation}
		\label{eqn-generalized_functional_invariant}
		[v_0,v_1] \mapsto [c_{ij}v_1^{i+j}\mathsf{t} : v_0^i(v_0+v_1)^{j}],
	\end{equation}
	for some constant $c_{i j} \in \mathbb{C}^\times$.  For a family $\pi \colon \mathcal{V} \to \mathbb{P}^1$ of elliptic CY varieties with Weierstrass models given by Equation~(\ref{eq-Weierstrass_model_family}) with complex $m$-dimensional fibers and a generalized functional invariant $(i, j, \alpha)$ such that 
	\begin{equation}
		\label{eqn-functional_invariant_conditions}
		0 \leq \mathrm{deg}_t(g_2) \leq \mathrm{min} \left ( \frac{4}{i},\frac{4\alpha}{j} \right ), \hspace{5mm} 0 \leq \mathrm{deg}_t(g_3) \leq \mathrm{min} \left ( \frac{6}{i},\frac{6\alpha}{j} \right ).
	\end{equation}
	Doran and Malmendier showed that a new family $\tilde{\pi} \colon \tilde{\mathcal{V}} \to \mathbb{P}^1$ can be constructed such that the general fiber $\tilde{V}_{\mathsf{t}} = \tilde{\pi}^{-1}(\mathsf{t})$ is a compact, complex $(m+1)$-manifold equipped with a Jacobian elliptic fibration over $\mathbb{P}^1 \times S$.  In the coordinate chart $\{ [ v_0 : v_1 ], (u_1,\dots, u_{m-1})\} \in \mathbb{P}^1 \times S$ the family of Weierstrass models $W_{\mathsf{t}}$, in the affine chart $\tilde{z}=1$, is given by 
	\begin{equation}
		\label{eqn-twisted_weierstrass}
		\begin{split}
			\tilde{y}^2  = 4\tilde{x}^3 & - g_2 \left ( \frac{c_{ij} \mathsf{t} v_1^{i+j}}{v_0^i(v_0 + v_1)^j},u \right ) v_0^4v_1^{4-4\alpha}(v_0+v_1)^{4\alpha} \tilde{x} \\
			& - g_3 \left ( \frac{c_{ij} \mathsf{t} v_1^{i+j}}{v_0^i(v_0 + v_1)^j},u \right ) v_0^6v_1^{6-6\alpha} (v_0+v_1)^{6\alpha}
		\end{split}
	\end{equation}
	with $c_{i j} = (-1)^i i^i j^j / (i+j)^{i+j}$ and 
	\begin{equation*}
		\tilde{x}=v_0^2 v_1^{2-2 \alpha}\left(v_0+v_1\right)^{2 \alpha} x\,, \qquad \tilde{y}=v_0^3 v_1^{3-3 \alpha}\left(v_0+v_1\right)^{3 \alpha} y.  
	\end{equation*}
	The new family is called the \emph{twisted family with generalized functional invariant $(i, j, \alpha)$ of $\pi : \mathcal{V} \to \mathbb{P}^1$}. It follows that conditions (\ref{eqn-functional_invariant_conditions}) guarantee that the twisted family is minimal and normal if the original family is. Moreover, they showed that the CY condition holds for the twisted family if it holds for the original.
	\par Additionally, a holomorphic $(m+1)$-form $\tilde{\Omega}_{\mathsf{t}} \in H^{m+1,0}(\tilde{V}_{\mathsf{t}})$ for the twisted family, in the affine chart $[v_0 \, : \, v_1] = [v \, : \, 1]$ and $\tilde{z} =1$, is given by
	\begin{equation}
		\label{eq-twisted_holomorphic_form}
		\tilde{\Omega}_{\mathsf{t}}=\frac{d v}{v(v+1)^\alpha} \wedge \Omega_{\mathsf{t}} \, ,
	\end{equation}
	where $\Omega_{\mathsf{t}}$ is the pullback of the holomorphic $m$-form in Equation (\ref{eq-holomorphic_form_WM}), pulled back via the map defining the generalized functional invariant in Equation (\ref{eqn-generalized_functional_invariant}). It follows that if $\omega(t)$ is a period integral for the family $\pi \colon \mathcal{V} \to \mathbb{P}^1$, a period integral $\tilde{\omega}(\mathsf{t})$ for the twisted family is as follows:
	\begin{equation}
		\label{eq-twisted_period}
		\tilde{\omega}(\mathsf{t})=\oint_C \frac{d v}{v(v+1)^\alpha} \, \omega\left(\frac{c_{i j} \mathsf{t}}{v^i(v+1)^j}\right) \, ,
	\end{equation}
	where $C$ is a suitably small non-contractible loop in the $v$-plane. Doran and Malmendier  showed from this integral representation that  the period integral for the twisted family with $\alpha = 1$ is the following Hadamard product  (\cite[Proposition 5.1]{MR4069107}):
	\begin{equation}
		\label{eq-Hadamard_product}
		\tilde{\omega}(\mathsf{t}) = (2\pi i) \, _{i+j-1}F_{i+j-2} \left(
		\begin{array}{ccc}
			\frac{1}{i+j} & \cdots & \frac{i+j-1}{i+j} \\
			\frac{1}{i} \cdots & \frac{i-1}{i} \, \frac{1}{j} & \cdots \frac{j-1}{j}
		\end{array} \Bigg{|} \mathsf{t}\right) \star \omega(\mathsf{t}) \, .
	\end{equation}
	A similar formula holds for $\alpha = \frac{1}{2}$, though we will not need it in this article. Here, $_pF_q$ is the generalized univariate hypergeometric function, and the Hadamard product $(f \star g)(t)$ of two power series $f(t) = \sum_{n \geq 0} a_n t^n$ and $g(t) = \sum_{n \geq 0} b_n t^n$ is given by  $(f \star g)(t) = \sum_{n \geq 0} a_n b_n t^n$. 
	
	\par We will be primarily interested in the twist construction for the generalized functional invariant $(i, j, \alpha)=(1,1,1)$. In this case, the relevant hypergeometric function that appears as the twist function in the Hadamard product is
	\begin{equation}
		\label{eq-1F0}
		_1F_0 \left(\frac{1}{2} \, \Bigg{|} \, 4t\right) = \frac{1}{\sqrt{1-4t}} = \sum_{k=0}^{\infty} \binom{2k}{k} \, t^k \, ,
	\end{equation}
	where we have chosen the argument $4t$ in order to produce an integral power series. This holomorphic function is annihilated by the differential operator 
	\begin{equation}
		\label{eq-1F0_op}
		\mathcal{L}_1 = (4t-1)\frac{d}{dt}+2 = t(2+4\theta)-\theta \, .
	\end{equation}
	These results prove that the Picard-Fuchs operator $\mathcal{L}_{n+1}$ of the twisted family $\tilde{\mathcal{V}} \to \mathbb{P}^1$ decomposes as the Hadamard product 
	\begin{equation}
		\mathcal{L}_{n+1} = \mathcal{L}_1 \star \mathcal{L}_n
	\end{equation}
	of the hypergeometric operator $\mathcal{L}_1$ and Picard-Fuchs operator $\mathcal{L}_n$ for the family $\mathcal{V} \to ~\mathbb{P}^1$. One also has that $\mathcal{L}_{n+1}$ is MUM at $t=0$ if $\mathcal{L}_n$ is. A similar statement holds for any generalized functional invariant $(i, j, \alpha)$.
	\newpage

	\subsection{Tables of modular data}
	\label{ss-modular_data}
	\par In this section we list tables of the relevant modular data and Picard-Fuchs operators under consideration in this article.
	\renewcommand{\arraystretch}{1.2}
	\begin{table}[!ht]
		\caption{\tiny Modular data for canonical elliptic fibrations over rational modular curves. $\eta(\tau)$ is the usual Dedekind $\eta$-function (\ref{eq-dedekind_eta}).  $\chi_{-3}(n),\chi_{-4}(n)$ are the quadratic Dirichlet characters $\chi_{-3}(n) = 0,1,-1 \; \; \mathrm{if} \; \; n \equiv 0,1,2 \; \mathrm{mod} \; 3$, $\chi_{-4}(n) = 0,1,0,-1 \; \; \mathrm{if} \; \; n \equiv 0,1,2,3 \; \mathrm{mod} \; 4$.  $\left(\frac{n}{5}\right)$, $\chi(n)$ are respectively the quadratic and quartic characters $\left(\frac{n}{5}\right) = (-1)^j, \; \chi(n) = i^j \; \; \mathrm{if} \; \; n \equiv 2^j \mod 5$. The Weierstrass data for $X_1(5)$ is the family that was utilized by Beukers \cite{MR0702189} in analyzing the irrationality of $\zeta(2) = \frac{\pi^2}{6}$. Note that the first two rows are actually orbifoldings of Beauville families.}\label{table-modular_elliptic_data}
		\scalebox{.8}{
			\begin{tabularx}{1.2\textwidth} { 
					>{\hsize=10mm\centering\arraybackslash}X 
					| >{\hsize=70mm\centering\arraybackslash}X |  >{\hsize=26mm\centering\arraybackslash}X |  >{\centering\arraybackslash}X  }
				$\Gamma \backslash \mathbb{H}^*$ & $g_2,g_3$ & $t(\tau)$ & $\omega_0(\tau)$,   $E_{3,m}(\tau)$   \\
				\hline 
				$X_0(3)$ & \tiny
				$\begin{aligned}
					g_2 &= 27 - 648 t \\
					g_3 &= 5832 t^2 - 972 t + 27
				\end{aligned}$	& $\frac{\eta(3 \tau)^{12}}{\eta(\tau)^{12}+27 \eta(3 \tau)^{12}}$  & \tiny  $\begin{aligned}
					\omega_0(\tau) &=1+6 \sum_{n=1}^{\infty} \chi_{-3}(n) \frac{ q^n}{1-q^n} \\
					E_{3,3}(\tau)&=1+9 \sum_{n=1}^{\infty}\chi_{-3}(n) \frac{ n^2q^n}{1-q^n}
				\end{aligned}$  \\
				\hline
				$X_0(4)$ & \tiny $\begin{aligned}
					g_2 &=\frac{1024}{3} t^2-\frac{1024}{3} t+\frac{64}{3} \\
					g_3 &= \frac{512}{27}(-1+8 t)\left(8 t^2+16 t-1\right)
				\end{aligned}$ & $ \frac{\eta(\tau)^8 \eta(4 \tau)^{16}}{\eta(2 \tau)^{24}} $ & \tiny $\begin{aligned}
					\omega_0(\tau) &=  1+4 \sum_{n=1}^{\infty} \chi_{-4}(n) \frac{ q^n}{1-q^n} \\
					E_{3,4}(\tau) &= 1+4 \sum_{n=1}^{\infty}\chi_{-4}(n) \frac{ n^2q^n}{1-q^n}
				\end{aligned}$ \\
				
				\hline
				$X_1(5)$ & \tiny $\begin{aligned}
					g_2 &= \frac{1}{12}\left(t^4 + 12t^3 + 14t^2 - 12t +1\right) \\
					g_3 &= -\frac{1}{216}\left(t^2 + 1\right)\left(t^4 + 18t^3  + 74t^2 - 18t +1\right) 
				\end{aligned}$ & \tiny $q \prod_{n=1}^{\infty}\left(1-q^n\right)^{5\left(\frac{n}{5}\right)}$ & \tiny $\begin{aligned}
					\omega_0(\tau) &=  1 + \frac{1}{2} \sum_{n=1}^{\infty} a_n\frac{q^n}{1-q^n} \\
					a_n &= (3-i)\chi(n) + (3+i)\overline{\chi(n)} \\
					E_{3,5}(\tau) &= 1+\frac{1}{2} \sum_{n =1}^{\infty}b_n \frac{n^2 q^n}{1-q^n} \\
					b_n &= (2-i) \chi(n)+(2+i) \overline{\chi(n)}
				\end{aligned}$ \\
				
				\hline
				$X_0(6)$ & \tiny $\begin{aligned}
					g_2 &= \frac{3}{4}(3 t+1)\left(3 t^3+75 t^2 -15 t+1\right)  \\
					g_3 &= \frac{1}{8}\left(3 t^2+6 t-1\right)  \left(9 t^4-540 t^3+30 t^2-12 t+1\right) 
				\end{aligned}$ & $\frac{\eta(\tau)^4 \eta(6 \tau)^8}{\eta(2 \tau)^8 \eta(3 \tau)^{4}}$ & \tiny $\begin{aligned} 
					\omega_0(\tau) &= \frac{\eta(2 \tau)^6 \eta(3 \tau)}{\eta(\tau)^3 \eta(6 \tau)^2} \\
					E_{3,6}(\tau) &= \frac{1}{9}\left(E_{3,3}(\tau)+8E_{3,3}(2 \tau)\right)
				\end{aligned}$ \\
				
				\hline 
		\end{tabularx}}
	\end{table}
	\renewcommand{\arraystretch}{1.2}
	\begin{table}[!ht]
		\caption{\tiny Modular Picard-Fuchs operators attached to modular elliptic surfaces and their associated quantum differential operators.}\label{table-modular_pf_elliptic}
		\scalebox{0.8}{
			\begin{tabularx}{.8\textwidth} { 
					>{\hsize=15mm\centering\arraybackslash}X 
					| >{\centering\arraybackslash}X  }
				$\Gamma \backslash \mathbb{H}^*$ & Picard-Fuchs operator $\mathcal{L}_{2,n}$, Quantum operator $\mathcal{D}_{3,n}$   \\
				\hline 
				\multirow{2}{4em}{$X_0(3)$} & $\theta^2 - 3t (3\theta + 1)(3\theta + 2)$, \\
				& $\theta^3-3 t(\theta+1)(3 \theta+1)(3 \theta+2)$ \\
				\hline
				\multirow{2}{4em}{$X_0(4)$} & $\theta^2 - 4t (2\theta + 1)^2$, \\
				& $\theta^3-4 t(\theta+1)(2 \theta+1)^2$\\
				\hline
				\multirow{2}{4em}{$X_1(5)$} & $\theta^2-t(11\theta^2+11\theta+3)-t^2(\theta+1)^2$, \\
				& $\theta^3-t(\theta+1)\left(11 \theta^2+11 \theta+3\right) -t^2(\theta+2)(\theta+1)^2$ \\
				\hline
				\multirow{2}{4em}{$X_0(6)$} & $\theta^2-t(10\theta^2+10\theta+3) +9t^2(\theta+1)^2$,  \\
				& $\theta^3-t(\theta+1)\left(10 \theta^2+10 \theta+3\right)+9 t^2(\theta+2)(\theta+1)^2$ \\
				\hline 
				%$X_0(8)$ & $\theta^2 - 4t(3\theta^2+3\theta+1)+32t^2(\theta+1)^2$ \\
				%\hline
				%$X_0(9)$ & $\theta^2 + 3t(3\theta^2+3\theta+1)+27t^2(\theta+1)^2$ \\
				%\hline
		\end{tabularx}}
	\end{table}
	\renewcommand{\arraystretch}{1.2}
	\begin{table}[!ht]
		\caption{\tiny Modular Picard-Fuchs operators and their associated quantum differential operators for rank-1 Fano threefolds of type $(N,1)$ that can be obtained naturally from the twist construction from the second order modular operators in Table \ref{table-modular_pf_elliptic}, for $N = 2,3,4,5$. The Ap\'ery family corresponds to $N=6$, and relevant twisted operator in $N = \widetilde{6}$. Hadamard products of the third order operators $\mathcal{L}_{3,N}$ are given in the second column.}\label{table-modular_pf}
		\scalebox{.8}{
			\begin{tabularx}{\textwidth} { 
					>{\hsize=5mm\arraybackslash}X 
					| >{\hsize=15mm\centering\arraybackslash}X 
					| >{\centering\arraybackslash}X  }
				$N$ & H. P. & Picard-Fuchs operator $\mathcal{L}_{3,N}$, Quantum Operator $\mathcal{D}_{4,N}$   \\
				\hline 
				\multirow{2}{4em}{2} & \multirow{2}{4em}{$\mathcal{L}_3 \star \mathcal{L}_{2,3}$} & $\theta^3 - 8t (2\theta + 1)(4\theta + 1)(4\theta + 3)$, \\
				& & $\theta^4-8 t(\theta+1)(2 \theta+1)(4 \theta+1)(4 \theta+3)$ \\
				\hline
				\multirow{2}{4em}{3} & \multirow{2}{4em}{$\mathcal{L}_1 \star \mathcal{L}_{2,3}$}  & $\theta^3 - 6t (2\theta + 1)(3\theta + 1)(3\theta + 2)$, \\
				& & $\theta^4-6 t(\theta+1)(2 \theta+1)(3 \theta+1)(3 \theta+2)$ \\
				\hline
				\multirow{2}{4em}{4} & \multirow{2}{4em}{$\mathcal{L}_1 \star \mathcal{L}_{2,4}$} & $\theta^3 - 8t (2\theta + 1)^3$, \\
				& & $\theta^4-8 t(\theta+1)(2 \theta+1)^3$  \\
				\hline
				\multirow{2}{4em}{5} & \multirow{2}{4em}{$\mathcal{L}_1 \star \mathcal{L}_{2,5}$} & $\theta^3 - 2t (2\theta + 1)(11\theta^2 + 11\theta + 3) - 4t^2(\theta+1)(2\theta + 1)(2\theta + 3)$, \\
				& & $\theta^4-2 t(\theta+1)(2 \theta+1)\left(11 \theta^2+11 \theta+3\right)-4 t^2(\theta+1)(\theta+2)(2 \theta+1)(2 \theta+3)$ \\
				\hline
				\multirow{2}{4em}{6} & \multirow{2}{4em}{\, $\times \times \times$} & $\theta^3 - t(2\theta+1)(17\theta^2+17\theta+5) + t^2(\theta+1)^3$,  \\
				& & $\theta^4 - t(\theta+1)(2\theta+1)(17\theta^2+17\theta+5)+ t^2(\theta+2)(\theta+1)^3$  \\
				\hline 
				\multirow{2}{4em}{$\widetilde{6}$} & \multirow{2}{4em}{$\mathcal{L}_1 \star \mathcal{L}_{2,6}$} & $\theta^3-2t(2\theta+1)(10\theta^2+10\theta+3)+36t^2(\theta+1)(2\theta+1)(2\theta+3)$, \\
				& & $\theta^4-2t(\theta+1)(2\theta+1)(10\theta^2+10\theta+3) +36t^2(\theta+1)(\theta+2)(2\theta+1)(2\theta+3)$\\
				\hline
		\end{tabularx}}
	\end{table}
	\renewcommand{\arraystretch}{1.2}
	\begin{table}[!ht]
		\caption{\tiny Modular data for mirror K3s associated rank-1 Fano threefolds of type $(N,1)$ that can be obtained from the twist construction from the modular elliptic surfaces in Table \ref{table-modular_elliptic_data}. Note that $N=6$ corresponds to the Ap\'ery family while $N=\widetilde{6}$ corresponds to the $\mathrm{M}_6$-polarized family, so as to be consistent with the Fano invariants. Here $G_2(\tau) = -\frac{1}{24}E_2(\tau)$ is the Hecke normalized Eisenstein series of weight-2 (\ref{eq-eisenstein_series_Ek}). While only a quasi-modular form, the combinations above in the period expressions are genuine modular forms of weight-2 for the corresponding modular group \cite{MR3462676}. Our normalizations for both differ slightly from loc. cit..}\label{table-K3_modular}
		\scalebox{.8}{
			\begin{tabularx}{\textwidth} { 
					>{\hsize=5mm\arraybackslash}X 
					| >{\hsize=52mm\arraybackslash}X 
					| >{\hsize=92mm\centering\arraybackslash}X  }
				$N$ & \; \; \; \; \; \; \; \; $t(\tau)$ & $\omega_0(\tau)$, $E_{4,N}(\tau)$   \\
				\hline 
				\multirow{2}{4em}{2} & \multirow{2}{0em}{$\frac{\eta(\tau)^{24} \eta(2 \tau)^{24}}{\eta(\tau)^{48}+128 \eta(\tau)^{24} \eta(2 \tau)^{24}+4096 \eta(2 \tau)^{48}}$} & \tiny $\omega_0(\tau) = 24G_2(\tau)- 2\cdot 24G_2(2\tau)$, \\
				& & \tiny $E_{4,2}(\tau)=-80 G_4(\tau)+2^2 \cdot 80 G_4(2 \tau)$ \\
				\hline
				\multirow{2}{4em}{3} & \;  \multirow{2}{4em}{$\frac{\eta(\tau)^{12} \eta(3 \tau)^{12}}{\eta(\tau)^{24}+54 \eta(\tau)^{12} \eta(3 \tau)^{12}+729(3 \tau)^{24}}$}  & \tiny $\omega_0(\tau) = 12G_2(\tau)- 3\cdot 12G_2(3\tau)$, \\
				& & \tiny $E_{4,3}(\tau)=-30 G_4(\tau)+ 3^2 \cdot 30  G_4(3 \tau)$ \\
				\hline
				\multirow{2}{4em}{4} & \; \; \; \; \; \multirow{2}{4em}{$\frac{\eta(\tau)^8 \eta(4 \tau)^8}{\left(\eta(\tau)^8+16 \eta(4 \tau)^8\right)^2}$} & \tiny $\omega_0(\tau) =\frac{\eta(2 \tau)^{20}}{\eta(\tau)^8 \eta(4 \tau)^8}$ , \\
				& & \tiny $E_{4,4}(\tau)=-16 G_4(\tau)+4^2 \cdot 16 G_4(4 \tau)$  \\
				\hline
				\multirow{2}{4em}{5} & \;  \multirow{2}{4em}{$\frac{\eta(\tau)^6 \eta(5 \tau)^6}{\eta(\tau)^{12}+22 \eta(\tau)^6 \eta(5 \tau)^6+125 \eta(5 \tau)^{12}}$} & \tiny $\omega_0(\tau) = 6G_2(\tau)-5\cdot6G_2(5\tau)$, \\
				& & \tiny $E_{4,5}(\tau)=-10 G_4(\tau)+  5^2 \cdot 10 G_4(5 \tau)$ \\
				\hline
				\multirow{2}{4em}{6} & \; \; \; \; \; \;  \multirow{2}{4em}{$\frac{\eta(\tau)^{12} \eta(6 \tau)^{12}}{\eta(2 \tau)^{12} \eta(3 \tau)^{12}}$} & \tiny $\omega_0(\tau) =\frac{\eta(2 \tau)^7 \eta(3 \tau)^7}{\eta(\tau)^5 \eta(6 \tau)^5}$,  \\
				& & \tiny $E_{4,6}(\tau)=9\left(-\frac{7}{9}G_4(\tau) + 2^2 \cdot \frac{1}{9}G_4(2\tau)- 3^2 \cdot \frac{1}{9}G_4(3\tau) + 6^2 \cdot \frac{7}{9}G_4(6\tau)\right)$  \\
				\hline 
				\multirow{2}{4em}{$\widetilde{6}$} & \; \; \multirow{2}{4em}{$\frac{\eta(\tau)^{4} \eta(2\tau)^{4} \eta(3 \tau)^4 \eta(6 \tau)^4}{\left(\eta(\tau)^4 \eta(2 \tau)^4+9 \eta(3 \tau)^4 \eta(6 \tau)^4\right)^2}$} & \tiny $\omega_0(\tau) = \frac{\eta(\tau)^{12} \eta(6 \tau)^{12}+\eta(2 \tau)^{12} \eta(3 \tau)^{12}}{\eta(\tau)^5 \eta(2 \tau)^5 \eta(3 \tau)^5 \eta(6 \tau)^5}$, \\
				& & \tiny $E_{4,\widetilde{6}}(\tau)=72\left(\frac{1}{9}G_4(\tau) - 2^2 \cdot \frac{1}{9}G_4(2\tau)+ 3^2 \cdot \frac{1}{9}G_4(3\tau) - 6^2 \cdot \frac{1}{9}G_4(6\tau)\right)$\\
				\hline
		\end{tabularx}}
	\end{table}

\end{document}